\documentclass{amsproc}

\usepackage{amsmath}
\usepackage{amsfonts}    
\usepackage{amssymb}
\usepackage{latexsym}
\usepackage{a4}
\usepackage{url}
\usepackage{graphicx}

\usepackage[x11names, rgb]{xcolor}
\usepackage{tikz}
\usetikzlibrary{snakes}
\usetikzlibrary{arrows}

\newcommand{\R}{\mathbb{R}}

\newcommand{\MT}{{\mathcal T}}

\newcommand{\MP}{{\mathcal P}}
\newcommand{\MQ}{{\mathcal Q}}

\newcommand{\MF}{{\mathcal F}}
\newcommand{\MH}{{\mathcal H}}

\renewcommand{\MR}{{\mathcal R}}

\newcommand{\SE}{\mathsf{E}}

\DeclareMathOperator{\Aut}{Aut}
\DeclareMathOperator{\Stab}{Stab}
\DeclareMathOperator{\Sym}{Sym}

\DeclareMathOperator{\Id}{Id}
\DeclareMathOperator{\conv}{conv}
\DeclareMathOperator{\aff}{aff}
\DeclareMathOperator{\cone}{cone}
\DeclareMathOperator{\gl}{GL}
\DeclareMathOperator{\gldr}{\gl_d(\R)}
\DeclareMathOperator{\push}{push}
\DeclareMathOperator{\pull}{pull}

\DeclareMathOperator{\argmin}{argmin}

\newtheorem{theorem}{Theorem}[section]
\newtheorem{definition}[theorem]{Definition}

\newtheorem{proposition}[theorem]{Proposition}
\newtheorem{example}[theorem]{Example}

\newcommand{\set}[1]{\{\,#1\,\}}
\newcommand{\one}{\mathbf{e}}

\newcommand{\ImName}{Incidence Decomposition Method}
\newcommand{\spname}{decomposition}
\newcommand{\SpName}{Decomposition}

\author{David Bremner}
\address{Faculty of Computer Science, UNB,
Box 4400, Fredericton, NB, E3B 5A3, Canada}
\email{bremner@unb.ca}
\author{Mathieu Dutour Sikiri\'c}
\address{Laboratory of Radiochemistry, Department of Physical Chemistry, Rudjer Boskovi\'c Institute, Bijenicka 54, 10000 Zagreb, Croatia}
\email{mathieu.dutour@ens.fr}
\author{Achill Sch\"urmann}
\address{Mathematics Department, University of Magdeburg, 39106 Magdeburg, Germany}
\email{achill@math.uni-magdeburg.de}

\thanks{The first author acknowledges the support of the Natural
  Sciences and Engineering Research Council of Canada and the
  Alexander von Humboldt Foundation.  The third author was supported
  by the Deutsche Forschungsgemeinschaft (DFG) under grant SCHU
  1503/4-1.  }

\title{Polyhedral representation conversion up to symmetries}

\begin{document}

\begin{abstract}
We give a short survey on computational techniques
which can be used to solve the representation conversion 
problem for polyhedra up to symmetries. 
We in particular discuss \spname{} methods, 
which reduce the problem to a number of lower dimensional subproblems. 
These methods have been successfully used by different authors in 
special contexts. 
Moreover, we sketch an incremental method, which is
a generalization of  Fourier--Motzkin elimination, and we give some
ideas how symmetry can be exploited using pivots.
\end{abstract}

\maketitle

\section{Introduction}

By the Farkas-Minkowski-Weyl Theorem 
a convex polyhedron in $\R^d$
has two representations.
It can either be described by a finite set of
linear inequalities (facets) or by a finite set of 
generators (vertices and rays).
Precise definitions are given in Section \ref{sec:basic_notations}.

One of the most fundamental problems in the theory of
polyhedra and its applications, such as 
Combinatorial Optimization or Computational Geometry,
is the conversion between two different descriptions.
Many algorithms for this representation conversion  have been proposed 
(see for example \cite{mr-1980}, \cite{abs-1997}, \cite{bfm-1998}, \cite{jz-2004}). 
For certain classes of polyhedra efficient methods are known,
but there is no approach known which efficiently 
solves the problem in general.
Programs like {\tt cdd} \cite{cdd}, {\tt lrs} \cite{lrs},
{\tt pd} \cite{pd}, {\tt porta} \cite{porta} 
(or {\tt polymake} \cite{polymake} either relying on some of the
others or using its own method closely related to \texttt{cdd})
allow conversion of the representation of a polyhedron. 
Since the programs are implementations
of quite different methods, their efficiency may vary tremendously on
a given example.

Many interesting polyhedra both pose difficulties for the standard
representation conversion approaches and have many symmetries that
could potentially be exploited.  In many applications it 
is sufficient (or at least necessary) to
obtain a list of inequalities or generators up to symmetries.  In the
present paper we give a brief survey of approaches that can be used
for the representation conversion problem up to symmetries.  We do not
discuss their asymptotic complexity (which, in the worst case, is not
encouraging), but rather refer to previous papers where the methods
have proven themselves on difficult instances that could not otherwise
be solved.  For the new approach discussed in Section
\ref{sec:sym-pivoting} we provide some experimental data ourselves.

The paper is organized as follows.
In Section \ref{sec:basic_notations} we give some 
basic notations and facts from the theory of convex polyhedra. 
In Section \ref{sec:symmetries} we consider different notions
of symmetries and describe how they
can practically be obtained using a graph automorphism computation.
In Section \ref{sec:orbits} we describe the group theoretical notions
used in the representation conversion methods discussed in the
remainder of the paper.
In Section \ref{sec:subpolytope} we consider
\spname{} methods which reduce the given problem to
a number of smaller problems. These approaches have
been used quite successfully by different authors.
In Section \ref{sec:cascade} we describe the 
incremental {\em Cascade algorithm} and
in Section \ref{sec:sym-pivoting} we show how 
symmetry can be exploited in a simplex pivot based algorithm.

\section{Convex polyhedra}

\label{sec:basic_notations}

In this section, we give a brief introduction to
some basic concepts and terminology of (convex) polyhedra. 
For  more details on polyhedra, 
we refer to the books \cite{ziegler-1998}, \cite{gruenbaum-2003}, \cite{schrijver}.

Given the vector space $\R^d$, denote by $(\R^{d})^{*}$ its dual
vector space, i.e.\ the vector space of linear functionals on $\R^{d}$.
A {\em convex polyhedron} $\MP\subseteq \R^d$ can be defined by a
finite set of linear inequalities
$$
\MP=\{ x\in\R^d : f_i(x)\geq b_i , i=1,\ldots,m\}
$$
with $f_i\in (\R^{d})^{*}$ and $b_i\in\R$ for $i=1,\ldots,m$.
If the number of inequalities $m$ in the description 
is minimum, we speak of a non-redundant description. 
The dimension $\dim \MP$ of $\MP$ is the dimension of the smallest
affine subspace containing it. 
Under the assumption that $\MP$ is full-dimensional, i.e.\ $\dim \MP =d$, 
every inequality $i$ of a non-redundant description defines 
a {\em facet} $\{x\in \MP : f_i(x)=b_i\}$ of $\MP$, 
which is a $(d-1)$-dimensional convex polyhedron contained
in the boundary of $\MP$.

By the Farkas-Minkowski-Weyl Theorem
(see e.g.\ \cite{schrijver}, Corollary 7.1a),
$\MP$ can also be described by a finite set of generators:
\begin{eqnarray*}
\MP & = & \conv\{v_1,\ldots,v_k\} + \cone\{v_{k+1},\dots,v_n\}  
\\
  & = & \{ \sum_{i=1}^n \lambda_i v_i : \lambda_i\geq 0, \sum_{i=1}^k \lambda_i=1\} 
\end{eqnarray*}
where $v_i\in\R^d$ for $i=1,\ldots,n$.
If the number of generators is minimum, the description is again called \emph{non-redundant}.
In the non-redundant case, the
generators $v_i$, $i=1,\dots,k$, are called {\em vertices} 
and $\R_{\geq 0}v_i$, $i=k+1,\dots,n$, are the {\em extreme rays} of $\MP$.
In case $\MP$ is bounded we have $n=k$ and we speak of 
a {\em convex polytope}.

The {\em representation conversion} from a minimal set of generators into
a minimal set of linear functionals (or vice versa) is called
the {\em dual description problem}.
By using homogeneous coordinates, 
the general inhomogeneous problem stated above can be reduced to the 
homogeneous one where $b_i=0$ and $k=0$.
For example, we embed $\MP\in\R^d$ in
the hyperplane $x_{d+1}=1$ in $\R^{d+1}$ and consider the closure of
its conic hull. 
The so-obtained polyhedron $\MP'=\cone\{v'_1,\ldots,v'_n\}$,
with $v'_i=(v_i,1)$ for $i=1,\ldots,k$ and $v'_i=(v_i,0)$ for $i=k+1,\ldots n$,
is referred to as {\em polyhedral cone}.

By duality, the problem of converting a description by homogeneous 
inequalities into a description by extreme rays is equivalent to the 
opposite conversion problem.
So for simplicity we assume from now on 
that $\MP\subseteq \R^d$ is a {\em polyhedral cone} given by a 
minimal (non-redundant) set of generators (extreme rays).
If
$$
\MP=\cone\{v_1,\dots,v_n\}
$$
we say that $\MP$ is generated by $v_1,\dots,v_n\in \R^d$.

We want to find a minimal set 
$\{f_1,\dots,f_m\}\subset (\R^{d})^{*}$ with
$$
\MP=\{ x\in\R^d : f_i(x)\geq 0 , i=1,\ldots,m\} 
.
$$
By choosing a suitable projection, it is possible to reduce the
problem further to the case where $\MP$ is full-dimensional and does
not contain any non-trivial linear subspaces.  For example, if
$v_1,\dots,v_n\in \R^d$ span a $k$-dimensional linear subspace, we may
just choose (project onto) $k$ independent coordinates.  The
appropriate projection (i.e.\  equations of the linearity space) can be
found efficiently via Gaussian elimination.  
In other words, without loss of generality we may assume  
the $v_i$ span $\R^d$ ($\MP$ is full-dimensional) 
and the linear inequalities $f_i$ span $(\R^{d})^{*}$
($\MP$ is pointed).

A {\em face} of $\MP$ is a set $\{x\in \MP : f(x)=0 \}$ where
$f$ is an element of 
$$
\MP^\ast 
=
\{
f\in(\R^{d})^{*} : f(x)\geq 0  \mbox{ for all } x\in \MP
\}
,
$$ 
the polyhedral cone {\em dual} to $\MP$. Note that $(\MP^\ast)^\ast = \MP$.
The faces of a pointed polyhedral cone are themselves pointed polyhedral
cones. We speak of a $k$-face, if its dimension is $k$. 
The faces form a (combinatorial) lattice ordered by inclusion, 
the  {\em face lattice} of $\MP$.
The rank of a face in the lattice is given by its dimension.
Each face is generated by a subset of the generators of $\MP$
and therefore it is uniquely identified by some subset of $\{1,\ldots,n\}$.
In particular, the $0$-dimensional face is identified with the
empty set $\emptyset$ and $\MP$ itself with the full index set $\{1,\ldots,n\}$.
All other faces of $\MP$ are identified with some strict, non-empty subset 
$F\subset \{1,\ldots,n\}$.

We write $F\lhd F'$ for two faces of $\MP$ with $F\subset F'$ and
$\dim F=\dim F'-1$.  Two $k$-faces of $\MP$ are said to be
\textit{adjacent}, if they contain  a common $(k-1)$-face and are
contained in a common  $(k+1)$-face.  In
particular, two extreme rays ($1$-faces) are adjacent, if they
generate a common $2$-face and two facets are adjacent, if they share
a common $(d-2)$-face (a {\em ridge}).  By the properties of a
lattice, for two faces $F_1$ and $F_2$, not necessarily of the same
dimension, there is always a unique largest common face $F=F_1\cap
F_2$ contained in them, and a unique smallest face $F'$, containing
both.  Any sub-lattice $[F:F']$ of the face lattice consisting of all
faces containing $F$ and contained in $F'$ is known to be isomorphic
to the face lattice of some pointed polyhedral cone of dimension $\dim
F'-\dim F$.  Note that this is a special feature of face lattices of
polyhedra, which is not true for general lattices.  In particular the
{\em diamond property} holds for face lattices: Every sub-lattice of
rank~$2$ carries the combinatorics of a $2$-dimensional polyhedral
cone, namely a face $F$ of rank~$k-1$, a face $F'$ of rank $k+1$ and
two faces $F_1$ and $F_2$ of rank $k$.  The rank~$2$ sub-lattices are
in one-to-one correspondence to pairs of adjacent faces.

A \emph{polyhedral complex} $\Delta$ is a set of polyhedral cones (the
\emph{cells} of $\Delta$) satisfying the following two properties:
\begin{enumerate}
\item[(a)] If $\MP \in \Delta$ then every face of $\MP$ is also in $\Delta$.
\item[(b)] For all $\MP_i,\MP_j \in \Delta$, $\MP_i \cap \MP_j$ is a 
  face of both $\MP_i$ and $\MP_j$.
\end{enumerate}
The facets of a polyhedral cone $\MP$ form a polyhedral complex, the
\emph{boundary complex} of $\MP$.  Polyhedral complex $\Delta'$ is a
\emph{subdivision} of polyhedral complex $\Delta$ if every cell of
$\Delta$ is the union of cells in $\Delta'$ and every cell of
$\Delta'$ is contained in some cell of $\Delta$. A subdivision is
called a \emph{triangulation} if every cell is a simplicial cone,
i.e.\ is the conic hull of exactly $d$-extreme rays.  By the
homogeneous embedding of polytopes discussed above, we may equally
discuss boundary complexes, subdivisions, and triangulations for
polytopes.  Two faces (cells)  of a polytopal complex have empty
intersection exactly when the corresponding cones intersect only at
the origin.

\section{Polyhedral Symmetries} 

\label{sec:symmetries}

\subsection{Groups acting on polyhedra}

In this paper we are especially interested in the case where some
non-trivial group acts on the given polyhedron $\MP$.  The {\em
  combinatorial automorphism group} of $\MP$, generated by $n$ extreme
rays, is the subgroup of $\Sym(n)$ of all permutations (acting on ray
indices $\{1,\ldots,n\}$) which preserve the whole face lattice of
$\MP$.  Thus the combinatorial automorphism group acts not only on the
generating rays, but also on the set of facets (inequalities), and
more generally on the sets of faces of a given dimension.  Moreover,
it preserves the inclusion relation among faces.  It is known that the
combinatorial automorphism group of $\MP$ can be computed from the
incidence relations between extreme rays and facets (see
\cite{kaibelschwartz}).

Given generators $v_1,\dots,v_n\in\R^d$ of a polyhedral cone $\MP$ 
and a subgroup $G$ of the combinatorial automorphism group of $\MP$,
we want to obtain the facets of $\MP$ up to symmetry, that is, 
one representative $F\subset \{1,\ldots,n\}$
from each orbit under the action of $G$.

If no group or only a small group is given, we usually still want
to exploit as much symmetries of the given polyhedron as possible 
in the representation conversion. That is, ideally we would
like to compute the full combinatorial automorphism group.
However, we do not know how to compute it without computing the facets,
which is precisely the problem we want to solve.
So we have to settle for a compromise and work with
a more restricted type of automorphism.
Note though, that one may obtain the full
combinatorial automorphism group of $\MP$, after
the facets have been computed

In many cases the combinatorial automorphism group (or some nontrivial
subgroup) may reflect geometric symmetries, for example if it has a
representation as a subgroup of $\gldr$ acting naturally on $\R^d$ and
the polyhedral cone $\MP$.  The group of all matrices in $\gldr$
preserving $\MP$ is called {\em linear automorphism group} of $\MP$.
Since a linear automorphism permutes the set of extreme rays
$\{\R_{\geq 0}v_1,\dots,\R_{\geq 0} v_n\}$, we naturally obtain a
representation as a permutation group $G\leq \Sym(n)$.  It is
important to note here that although the induced permutation group is
finite, the linear automorphism group of $\MP$ is not necessarily so,
and it can be quite awkward. Think for example of the quadrant in
$\R^2$ generated by the two non-negative coordinate axes.  The induced
permutation group is $\Sym(2)$, its linear automorphism group however
is
$$
\left\{
\begin{pmatrix}
a & 0\\
0 & b\\
\end{pmatrix}
,
\begin{pmatrix}
0 & c\\
d & 0\\
\end{pmatrix}
:
a,b,c,d\in\R_{>0}
\right\}
.
$$

As explained in Section \ref{sec:basic_notations} we may limit our 
discussion to the special case of a full-dimensional, pointed polyhedral cone $\MP$
generated by $v_1 , \ldots, v_n \in\R^d$.
The elements $A\in\gldr$ of the linear automorphism group of $\MP$
satisfy $Av_i = \lambda_i v_{\sigma(i)}$, where $\sigma\in\Sym(n)$
is an induced permutation and $\lambda_i>0$.
Note, in case $\MP$ is the homogenization of a 
$(d-1)$-dimensional polyhedron $\MP'$, 
the linear automorphism group of $\MP$ corresponds to the 
so called \textit{projective automorphism group} of $\MP'$,
that is, the set of all projective maps preserving $\MP'$.
We are not aware of any practical algorithm to compute linear (or
projective) automorphism groups or to decide if two polyhedral cones
are linear (or projective) isomorphic.

\subsection{Restricted isomorphisms}
\label{sec:restrict}

A {\em restricted isomorphism} of two full-di\-men\-sio\-nal vector families
$V=\{v_1,\dots,v_n\}$ and $V'=\{v'_1,\dots,v'_n\}$ in $\R^d$ is given
by a matrix $A\in \gldr$ such that there exists a permutation $\sigma$
satisfying $Av_i=v'_{\sigma(i)}$ for $i=1,\dots, n$.  
A {\em restricted automorphism} of a vector family is a restricted 
isomorphism of $V$ with itself.

We speak of a restricted isomorphism between two polyhedral cones
generated by the two vector families $V$ and $V'$, if it is a
restricted isomorphism between $V$ and $V'$.
Note though that the definition of
restricted isomorphisms for polyhedral cones
depends strongly on the choice of generators.
Since such generators are only unique up to any positive factor,
the choice of generators is very crucial.
In practice the situation is usually not so bad 
as there is often a natural choice for the generators.

The big advantage of restricted isomorphisms is
that we can compute them by obtaining the graph isomorphisms 
for an edge colored graph.
Given a vector family $V=\{v_1,\dots, v_n\}\subset \R^d$ spanning $\R^d$
we define the positive definite matrix
\begin{equation}   \label{eqn:Q-def}
Q=\sum_{i=1}^{n} v_i v_i^t
.
\end{equation}
Let the graph $G(V)$ be the complete graph with vertices $v_i$
and edge colours $c_{ij}=v_i^{t} Q^{-1} v_j$.
Then the following holds:

\begin{proposition}  \label{prop:isomorphism-equivalence}
Let $V, V'\subset\R^d$ be two finite vector families.
Then every isomorphism of the edge coloured graphs
$G(V)$ and $G(V')$ yields a restricted isomorphism 
of $V$ and $V'$ and vice versa.
\end{proposition}

In practice the popular and very nice program {\tt nauty} \cite{nauty}
by McKay can be used to check for graph isomorphisms or compute the
group of automorphisms.  Note however, that the current version only
takes vertex coloured graphs as input.  We therefore have to transform
our edge coloured graph to a somewhat larger vertex coloured graph
that preserves the automorphism group (see for example \cite{nauty-manual},~p.25).

\begin{proof}[Proof of \ref{prop:isomorphism-equivalence}]
Let $Q$ be the matrix \eqref{eqn:Q-def} obtained from $V$.
Denote by $R$ the unique square root of $Q^{-1}$,
that is, the positive definite $d\times d$ matrix $R$ with $Q^{-1}=R^2$.
Let $w_i=Rv_i$ for $i=1,\ldots,n$. Then the edge colours $c_{ij}$
of the graph $G(V)$ are exactly the inner products 
$w_i^t w_j$ of the transformed vectors $w_i$ and $w_j$.
In the same way we obtain $Q'$, $R'$ and $w'_i$ for $V'$.

Now let $A\in\gldr$ be a restricted isomorphism for $V$ and $V'$
with associated permutation $\sigma\in\Sym(n)$.
Then we have $AQA^t = Q'$. 
Moreover, the matrix
$T=R'AR^{-1}$ is orthogonal
and satisfies $Tw_i=w'_{\sigma(i)}$.
This implies $c_{ij}=c'_{\sigma(i)\sigma(j)}$, i.e.\ the restricted
isomorphism of the vector families $V$ and $V'$ corresponds
to an isomorphism between the edge coloured graphs $G(V)$ and $G(V')$.

Suppose on the other hand $\sigma\leq\Sym(n)$ be an isomorphism
between the edge coloured graphs $G(V)$ and $G(V')$.
By reordering elements, we may simply assume $\sigma=\Id$.
Instead of looking for a solution of $Av_i=v'_i$ 
we consider the equivalent equations $Tw_i=w'_i$.
Since the $v_i$ generate $\R^d$, we find a basis
$(v_{i_1}, \dots, v_{i_d})$ of $\R^d$.
If $P$, respectively $P'$ is the $d\times d$ matrix formed by $(w_{i_k})$,
respectively $(w'_{i_k})$, then the matrix equation $(c_{ij})=(c'_{ij})$ takes the
form $P^tP=P'^tP'$. So, the matrix $T=P'P^{-1}$ is orthogonal
and we have for any $k=1,\dots,d$ and $j=1,\dots, n$:
\begin{equation*}
w'^t_{i_k}  Tw_j
=
(Tw_{i_k})^t Tw_j
=
w^t_{i_k} w_j
=
w'^t_{i_k} w'_j
.
\end{equation*}
This yields $w'^t_i (Tw_j-w'_j) = 0$ and since the $w'_i$
form a basis of $\R^d$ the relation $Tw_j=w'_j$ are implied.
Thus we obtain a restricted isomorphism between $V$ and~$V'$.
\end{proof}

\section{Orbits of faces}

\label{sec:orbits}

\subsection{Dealing with orbits}

In order to generate facets (or more generally faces)
of a polyhedron up to symmetries, it is necessary to deal 
with orbits of faces. In this section
we briefly indicate what the basic tasks are that we have to 
accomplish and how these can be approached.

As before, we assume that a polyhedral cone $\MP$ is given by a set of
generators $v_1,\dots,v_n\in\R^d$ which define the extreme rays of
$\MP$. Each face is represented by a subset of $\{1,\dots,n\}$ which
corresponds to the indices of generators incident to the face.  We
assume $G\leq \Sym(n)$ is some subgroup of the combinatorial
automorphism group of $\MP$, hence a permutation group acting not only
on the set of indices $\{1,\dots,n\}$ (respectively rays), but also on
the whole face lattice of $\MP$.  In particular, the dimension of a
face and inclusion between faces are preserved by every group element.

In the most general form, the problem we have to solve is the
following: Given two subgroups $G_1$ and $G_2$ of $\Sym(n)$
and a list $L_1$ of $G_1$-inequivalent faces, we need to
be able to obtain a list $L_2$ of $G_2$-inequivalent faces.

For example, if $G_1=\{\Id \}$ is trivial, the list $L_1=\{F_1,\dots,F_k\}$ 
simply is a list of pairwise unequal faces and we may obtain $L_2$ by
testing for $G_2$-equivalence: Starting with $L_2=\{F_1\}$ and
then subsequently adding $F_i$ for $i=2,\dots,n$ to $L_2$, 
if it is not $G_2$-equivalent to any element in $L_2$.
Clearly, this {\em orbit fusion} can be applied whenever
$G_1 \leq G_2$. 

In case $G_2 < G_1$ it is necessary to ``break some symmetry'' 
and {\em split (factorize) orbits}.
This can be done with the {\em double coset decomposition},
described in Section \ref{sec:canonical-rep}.
Note, if neither $G_1 \leq G_2$ nor $G_2 < G_1$, 
one could in principle convert the $G_1$-orbits 
into $G_2$-orbits in two steps,
by either using the intersection group
$G=G_1\cap G_2$ or the group $G=\langle G_1, G_2\rangle$
generated by $G_1$ and $G_2$ and obtaining 
$G$-orbits in an intermediate step.

 Orbit fusing and splitting is a typical and essential 
task when generating discrete structures up to isomorphism 
(see for example \cite{glm-1997}, \cite{kerber}, \cite{brinkmann}, \cite{ko-2006}).

\subsection{Fusing orbits, equivalence tests and canonical representatives} 

\label{sec:canonical-rep}

A common task is to decide whether or not two faces are $G$-equivalent, 
that is, whether or not the corresponding subsets 
of $\{1,\dots,n\}$ lie in the same orbit under the action of $G$.

The dimension of a face (rank in the face lattice) and its cardinality
are quickly testable invariants.  Other invariants can easily be
obtained, for example by looking at the action of $G$ on pairs,
triples, or other $k$-tuples of indices (generators), respectively on
lower dimensional faces.  The number of elements from each such orbit
included in a face is a $G$-invariant.  Unfortunately, there is no
clear rule how many such sets have to be chosen, so one has to rely on
heuristics.

If we only consider restricted isomorphisms, metric invariants
can be used as well, for example  
the set of pairwise inner products
between generators, discussed further in Section~\ref{sec:metric-invariants}.

When all invariants are satisfied, then group computations must be done.
For small groups it is possible to simply generate the whole orbit of
a face. If the size of available memory is a problem, we may just keep
a {\em canonical representative} for the orbit, 
for example the lexicographical minimum 
(when viewed as a subset of $\{1,\dots,n\}$). 
These can be found for example by a backtrack method 
(see Section~\ref{sec:implementing_orbits}).

\subsection{Breaking symmetry by splitting orbits using double cosets}

\label{sec:double-cosets}

Assume $G_2<G_1\leq \Sym(n)$ and that a list of $G_1$-orbits is given.
In order to obtain a list of $G_2$-orbits, we may split each
orbit $G_1 F$, with a representative face $F$, by the well known
{\em double coset decomposition} (cf.\ for example in \cite{brinkmann} and \cite{kerber}):
The group $G_1$ can be decomposed into {\em double cosets}
$$
G_1=\bigcup_{i=1}^r G_2 g_i\Stab(G_1,F)
$$
where $g_1, \dots, g_r$ are elements of $G_1$.
Then the orbit $G_1 F$ is decomposed into 
$$
G_1 F=\bigcup_{i=1}^r G_2 g_i F,
$$
hence into $r$ orbits $G_2 F_i$ with $F_i=g_iF$.

Note that the more straightforward algorithm of generating the full orbit
or of computing a decomposition of $G_1$ into cosets $G_2 g_i$ 
is slower and/or requires more memory (see \cite{kerber}).

\subsection{Data structures and implementation issues} 

\label{sec:implementing_orbits}

The fundamental data structures used to work with permutation groups
are {\em bases and strong generating sets} (BSGS) 
(see for example \cite{seress-2003}, \cite{heo-2005}).
Based on them a backtrack search on cosets of a point stabilizer
chain can be used to obtain canonical representatives, 
to decide on (non-)equivalence as well as to obtain stabilizers 
of subsets of $\{1,\dots,n\}$ (faces). 
For details we refer to \cite{seress-2003} and \cite{ko-2006}.
An elaborate version is the {\em partition backtrack} introduced
by Leon  \cite{leon-1991} (cf.\ \cite{leon-1997}, \cite{seress-2003}).
These methods are known to work quite well in practice, although
from a complexity point of view the problems are known to be difficult.
That is, there is no (worst-case) polynomial time algorithm known 
to solve these problems.
Even worse, it is somewhat unlikely that there exist polynomial time algorithms, 
since the {\em graph isomorphism problem}, not known to be in $\mathsf{P}$, 
is reducible to them in polynomial time (cf.\ \cite{luks-1993}).

The computer algebra system GAP \cite{GAP} provides functions
for generation of full orbits ({\sc Orbit}) 
stabilizer computations ({\sc Stabilizer}) 
and equivalence tests ({\sc RepresentativeAction}).

In the case where we have special knowledge of the group $G$ or its representation, 
it might be much easier to obtain canonical representatives
or to compute stabilizers of faces. 
For example, we may have a situation where the 
symmetric group $\Sym(n)$ acts on $n$ elements
(see \cite{anzin} and \cite{robbins}).
Another example where the action of the symmetric group is used
is described in \cite{di-2007} 
(see \cite{di-2006} for corresponding computer code).
Typically, polyhedra
arising in Combinatorial Optimization 
are convex hulls of $(0/1)$-vectors, 
where each coordinate (variable) represents an edge of a complete 
directed or undirected graph with $n$ vertices on which
$\Sym(n)$ acts.
In \cite{CR02} a method for obtaining
canonical representatives in this situation is described.

\section{\SpName{} methods}

\label{sec:subpolytope}

In this section we describe two methods which reduce 
the facet generation problem under symmetries
to a number of smaller instances of 
facet generation problems. 
In contrast to the original problem, 
solving the smaller problems (for sub-cones) may 
be feasible for available software, such as 
{\tt cdd} or {\tt lrs}.
These techniques have been proven to be successful
in practice, in cases where standard approaches failed.

We mainly distinguish two approaches, the {\em Incidence Decomposition Method} 
(see Section \ref{sec:incident}) 
and the {\em Adjacency Decomposition Method} 
(see Section \ref{sec:adjacency}).
Both methods are reasonably natural and have been used separately 
by different authors.

\subsection{\ImName}

\label{sec:incident}

The Incidence Decomposition Me\-thod reduces the problem of facet generation to
a number of smaller problems, in which we generate facets that are incident
to some extreme rays. As before let $\MP$ be a polyhedral cone in $\R^d$,
generated by $v_1,\dots, v_n$.
Assume  $G\leq \Sym(n)$ is some 
permutation group acting on the face lattice of $\MP$. 
The set of extreme rays (indices) falls into orbits 
under the action of $G$.
For each orbit we consider a representative $r_i$ (index $i$) and
generate a list of $G$-inequivalent facets of $\MP$ incident to it. 
Then, in a post-processing step the lists of
facets obtained in this way are merged to a list of 
$G$-inequivalent facets of $\MP$ (see Section \ref{sec:orbits}). 
Since every $G$-orbit of facets of $\MP$ contains a facet
which is incident to one of the chosen representatives, 
the resulting list is complete.

\begin{flushleft}
\smallskip
\textbf{Input:} $n$ extreme rays of a polyhedral cone $\MP$ and a group $G\leq\Sym(n)$ acting on $\MP$'s face lattice.\\
\textbf{Output:} complete set~$\MF$ of $G$-inequivalent facets of $\MP$.\\
\smallskip

$\MF \leftarrow \emptyset$.\\ 
$\MR \leftarrow \mbox{complete set of $G$-inequivalent extreme rays of $\MP$}$.\\
\textbf{for} $r \in \MR$ \textbf{do}\\
\hspace{2ex} $\MF_r \leftarrow \;\mbox{facets of $\MP$ incident to $r$}$.\\
\hspace{2ex} \textbf{for} $F \in \MF_r$ \textbf{do}\\
\hspace{2ex} \hspace{2ex} \textbf{if} $F$ is $G$-inequivalent to facets in $\MF$ \textbf{then}\\
\hspace{2ex} \hspace{2ex}\hspace{2ex} $\MF \leftarrow \MF \cup \{F\}$.\\
\hspace{2ex} \hspace{2ex} \textbf{end if}\\
\hspace{2ex} \textbf{end for}\\
\textbf{end for}\\
\end{flushleft}

The main computational gain comes from the following:
when we compute the facets that are incident to a given ray $r_i$,
we may not have to consider all $n$ extreme rays of $\MP$, because some of the rays
may not be incident to facets which are incident to $r_i$.
These are exactly the extreme rays $r_j=\R_{\geq 0}v_j$ which give redundant
inequalities of the polyhedron
$$
\MP_i^\ast
:=
\{
f\in (\R^d)^{\ast} 
:
f(v_j)\geq 0, j=1,\dots,n \;\mbox{and}\;
f(v_i)=0
\}
.
$$ 
So, for each of the $n-1$ rays with index in $\{1,\dots,n\}\setminus\{i\}$
we may solve a linear program in $(\R^d)^{\ast}$ to decide redundancy.
As an outcome we obtain a list of $n'<n$ rays needed in the definition of $\MP_i^\ast$.
The smaller $n'$, the bigger the computational gain.
Note, that $\MP_i^\ast$ is of dimension $d-1$. Its dual contains 
the line $\{\lambda v_i : \lambda\in\R \}$ and we may project it 
(or at least think of it as projected) down along this line to a $d-1$ dimensional
polyhedral cone.

So the problem of enumerating facets incident to a given ray is a
facet enumeration problem, but in one lower dimension and with fewer
extreme rays. Some of the lower dimensional subproblems may still be
too difficult; in this case we may apply the method recursively.  We
come back to this in Section~\ref{sec:recursion}.

The \ImName{} has been used in \cite{GrMet} for finding the vertices
of the metric polytope $\mbox{MET}_7$.  The method was also introduced
in \cite{CR} along with the Adjacency Decomposition Method but it was
found to be less competitive for their application. The \ImName{} was
also discussed by Fukuda and Prodon~\cite{fp96}.

\subsection{Adjacency Decomposition Method}

\label{sec:adjacency}

As with the \ImName{}, the Adjacency Decomposition Method 
is a reasonably natural method for computing the facets of a 
polytope up to symmetries. 
So it is no wonder that the method was discovered several
times, for example, in \cite{jaquet} as ``algorithm de l'explorateur'',
in \cite{CR} as ``adjacency decomposition method'' and
in \cite{DFPS} as ``subpolytope algorithm''.
Other example of applications are in \cite{rigid}, \cite{mining} or \cite{quasi2}.
The adjacency decomposition scheme is also implicit in polyhedral 
decomposition schemes, i.e.\ when a space is decomposed as an union 
of polyhedral cones. Examples are in \cite{sturmfels}, 
\cite{voronoi1} and \cite{voronoi2}
(cf.\ \cite{dsv-2006b} and \cite{dsv-2006a}).

Rather than focusing on the incidence of (orbits of) facets to extreme
rays as the \ImName{} does, the Adjacency Decomposition Method focusses
on the incidence of facets with other facets.
Starting from a (set of) initial $G$-inequivalent facets, it
traverses the adjacency graph of facet orbits. The initial facet(s) may
be obtained by suitable linear programs, or in strongly polynomial
time using the methods described in~\cite{bfm-1998}.

\begin{flushleft}
\smallskip
\textbf{Input:} $n$ extreme rays of a polyhedral cone $\MP$ and a group $G\leq\Sym(n)$ acting on $\MP$'s face lattice.\\
\textbf{Output:} complete set~$\MF$ of $G$-inequivalent facets of $\MP$.\\
\smallskip

$\MT \leftarrow \{F\}$ with $F$ a facet of $\MP$.\\
$\MF \leftarrow \emptyset$.\\
\textbf{while} there is an $F \in \MT$ \textbf{do}\\
\hspace{2ex} $\MF \leftarrow \MF \cup \{F\}$.\\
\hspace{2ex} $\MT \leftarrow \MT \setminus \{F\}$.\\
\hspace{2ex} $\MH\leftarrow \mbox{facets of $F$}$.\\
\hspace{2ex} \textbf{for} $H \in \MH$ \textbf{do}\\
\hspace{2ex} \hspace{2ex} $F' \leftarrow$ facet of $\MP$ adjacent to $F$ along $H$.\\
\hspace{2ex} \hspace{2ex} \textbf{if} $F'$ is $G$-inequivalent to all facets in $\MF\cup \MT$ \textbf{then}\\
\hspace{2ex} \hspace{2ex}\hspace{2ex} $\MT \leftarrow \MT \cup \{F'\}$.\\
\hspace{2ex} \hspace{2ex} \textbf{end if}\\
\hspace{2ex} \textbf{end for}\\
\textbf{end while}\\
\end{flushleft}

The facet $F'$ of $\MP$ with $F\cap F'=H$ for a given ridge $H\lhd F$
can be found by a \emph{gift-wrapping step} (cf.~\cite{ck-acp-70}, \cite{s-fchff-85}).  
Let $v_1,\dots,v_n$ be the
generators of $\MP$'s extreme rays.  The defining inequality $f\in
(\R^d)^\ast$ of the facet $F'$ should satisfy $f(v_i)=0$ for all
generators $v_i\in H$.  The vector space of such functions has
dimension $2$.  Let us select a basis $\{f_1, f_2\}$ of it.  If
$f=\alpha_1 f_1+\alpha_2 f_2$ is the defining inequality of $F$, $H$
or $F'$, then $f(v_i)\geq 0$ for all $i\in\{1,\dots,n\}$.  This
translates into a set of linear inequalities on $\alpha_1, \alpha_2$
defining a $2$-dimensional pointed polyhedral cone.  One finds easily
its two generators $(\alpha^{i}_1, \alpha^{i}_2)_{1\leq i\leq 2}$.
The corresponding inequalities $f_i=\alpha^i_1 f_1+\alpha^i_2 f_2\in
(\R^d)^\ast$ define the two adjacent facets $F$ and $F'$ of $\MP$.
In the special case where no $d+1$ extreme rays lie on a
hyperplane, the gift wrapping step corresponds to a simplex
pivot (see \cite{chvatal}).

\medskip

A nice feature of the Adjacency Decomposition Method is that the
adjacencies of the most symmetric facets, which are usually the most
difficult to treat, may not have to be considered. This is due to the
well known \emph{Balinksi's Theoreom} (see e.g. \cite{ziegler-1998},
Theorem 3.14):

\begin{theorem}\label{Balinski}(\cite{balinski})
Let $\MP$ be a $d$-dimensional, pointed polyhedral cone. Let $G$ be the
undirected graph whose vertices are the facets of $\MP$ and whose
edges are the ridges of $\MP$. Two vertices $F_1, F_2$
are connected by an edge $E$ if $F_1\cap F_2=E$.
Then, the graph $G$ is $(d - 1)$-connected, i.e.\ removal of any
$d - 2$ vertices leaves it connected.
\end{theorem}

Using Balinski's Theorem, we know that if the number of facets in
unfinished orbits (i.e.\ those whose neighbors are not known) is less than $d-1$, then any $G$-inequivalent facet
in the unfinished set must be adjacent to some already completely
treated facet. 
 This simple criterion has proven~\cite{dsv-2006b}, \cite{dsv-2007c} to be extremely
 useful in dealing with examples arising in the geometry of numbers.

\subsection{Recursion}

\label{sec:recursion}

Both, the \ImName{} and the Adjacency Decomposition Method
reduce the facet enumeration problem for a $d$-dimensional polyhedral cone $\MP$  
to a number of facet enumeration problems for cones in dimension $d-1$.
These lower dimensional problems may be too difficult to treat
with a standard method as well and therefore we might 
apply the Incidence or Adjacency Decomposition Method to these
lower dimensional problems recursively. So we may speak of the
{\em Recursive \ImName{}}, the {\em Recursive Adjacency Decomposition Method}
and the {\em Recursive Decomposition Method}, if a mixture of 
both is applied.

The Recursive \ImName{} has successfully been used  
for computing the vertices of the metric polytope 
$\mbox{MET}_8$ in \cite{dfmv}, whereas the
Recursive Adjacency Decomposition Method 
has successfully been used
in \cite{jaquet} and \cite{dsv-2006b}.
To the best of our knowledge, and to our surprise, 
a recursive mixture of both has not been used so far.

The crucial steps for the recursion are ``$\MF_r \leftarrow
\;\mbox{facets of $\MP$ incident to $r$}$'' (where $r$ is some extreme
ray) for the Incidence Method (see Section \ref{sec:incident}) and
``$\MF \leftarrow \mbox{facets of $F$}$'' (where $F$ is some facet)
for the Adjacency Decomposition Method (see Section
\ref{sec:adjacency}).  In both cases the problem is to obtain a list
of facets for a $(d-1)$-dimensional polyhedral cone $F$, whose extreme
rays are given. 
Again we can exploit symmetries.
As a result of a call of the Incidence or Adjacency Decomposition Method
for $F$, we obtain a list of $G_F$-inequivalent facets of $F$ (ridges of $\MP$),
where $G_F$ is some group acting on the face lattice of $F$
which we have to provide.
In a post processing step we then have to obtain a list of 
$G_{\MP}$-inequivalent (respectively $\Stab(G_{\MP},F)$-inequivalent)
facets of $F$ (ridges of $\MP$) out of it.

Note though,
that the groups $G_F$ and $G_\MP$ may be unrelated, that is,
the elements of $G_F$ do not have to be elements of $G_\MP$ and vice versa.
Think of examples where $\MP$ is a polyhedral cone without any symmetries,
but with a facet having some symmetries.
Vice versa, since we did not assume that $G_F$ is the full (combinatorial) symmetry group, 
not even the stabilizer $\Stab(G_{\MP}, F)$ of $F$ in $G_{\MP}$ has
to be a subgroup of $G_F$. 
In the latter case we may simply replace $G_F$ by the group generated
by $\Stab(G_{\MP}, F)$ and $G_F$, so that we assume $\Stab(G_{\MP}, F)\leq G_F$.
In this way we possibly enlarge the group and speed up the computations. 
Moreover, we are able to use the double coset decomposition 
(see Section \ref{sec:double-cosets}) to split each $G_F$-orbit
$G_F F'$ of facets of $F$ into a finite number of $\Stab(G_{\MP},F)$-orbits 
$\Stab(G_{\MP}, F)g_i F'$, where the $g_i\in G_F$ represent the double cosets.

Using the described \spname{} methods recursively it might happen that
we compute the facets of some sub-cones several times.  For example if
a face $F$ satisfies $F\lhd F_1\lhd \MP$ for a facet $F_1$ of a
polyhedral cone $\MP$, then there is exactly one other facet $F_2$
such that $F\lhd F_2\lhd \MP$. Hence, if we apply the Adjacency
Decomposition Method to $F_1$ and $F_2$, then we have to compute the
dual description of $F$ two times.  Clearly, the number of such
repetitions increases as the recursion depth increases. Moreover,
equivalent sub-cones may occur in different parts of the face lattice.
(Recall that any sub-lattice of the face lattice is the face lattice
of some polyhedral cone.)  To handle this, we propose to use a
\textit{banking system}.  That is, given a difficult sub-cone $\MP$ of
which the facets have been computed up to symmetries, we store with
$\MP$ (represented by generators) its group of restricted automorphisms
and a representative for each facet orbit.

What about the recursion depth needed to practically solve a given
polyhedral conversion problem? If we split the problem into
subproblems with either decomposition method then some of the
subproblems may be easy to solve and may not require a recursive
treatment, others might be impossible to treat without such.  The
implementation~\cite{Dut07} allows one to provide a heuristic function
to choose whether to recursively apply the Adjacency Decomposition
Method, but in choosing the right level of recursion requires some
trial and error.  From our computational experience with this code
(for example in \cite{rigid}, \cite{mining}, \cite{quasi2}, \cite{dsv-2006a}, \cite{dsv-2006b}), 
the {\em incidence number} of a face, that is, the number of
extreme rays contained in it, gives a good measure for how difficult
the representation conversion is for it.  We therefore propose to
treat the faces with low incidence numbers first and as much as
possible without recursion.  For subproblems to be potentially solved
with pivoting methods (possibly up to symmetry as below), the ``probing
feature'' of {\tt lrs} can be used to test for expected difficulties
of a subproblem.

Finally, let us remark, that 
parallelization is a possible way to speed up the (recursive) 
Decomposition Methods as well (see \cite{CR02} and \cite{di-2006}).

\section{An incremental method}

\label{sec:cascade}

We briefly sketch here the {\em Cascade algorithm} by Jaquet \cite{jaquet},
which is a symmetry exploiting version of {\em Fourier--Motzkin elimination}.
For this observe first that a $d$-dimensional polyhedral cone $\MP$ 
generated by $v_1 , \ldots, v_n \in\R^d$
can be obtained as the projection of an $n$-dimensional polyhedral
cone with $n$ linear independent generators $v'_1,\ldots,v'_n$.
Without loss of generality, we may assume that
$(v_1, \dots, v_d)$ is a basis of $\R^d$.
We set $v'_i=(v_i, 0^{n-d})$ for $i\leq d$ and $v'_i=(v_i, 0^{n-d})+e_i$ for $i > d$.
Let $p_i(x_1,\dots, x_n)=(x_1, \dots, x_i, 0,\dots,0)$ be the
orthogonal projection of $\R^n$ onto $\R^i\times \{0\}^{n-i}$.
Then for the cone $\MP'=\cone \{ v'_1, \ldots, v'_n \}$
we have $\MP=p_{d}(\MP')$. The cone $\MP'$ is simplicial, since
the family $(v'_i)_{1\leq i\leq n}$ is a basis of $\R^n$.

Take the $n$-dimensional polyhedral cone $\MP'$ and the projection $p_{n-1}$
of $\R^n$ on a $(n-1)$-dimensional hyperplane. The facets of the
projection $p_{n-1}(\MP')$ are either projections of facets of $\MP'$
or projections of intersections of facets of $\MP'$.
In the Fourier--Motzkin elimination we first compute the facets of $\MP'$
and then successively obtain the facets of  
$p_{n-1}(\MP'),\ldots, p_{d}(\MP')=\MP$ from one of these two cases.
Note, since $\MP'$ is a simplex, the facets of $\MP'$ are simply given 
by the $(n-1)$-subsets of extreme rays of $\MP'$.

In the Cascade algorithm we also consider symmetries in each step:
If $G$ is a group of symmetries of the polyhedral cone $\MP$, 
then the induced symmetry group on $p_i(\MP')$ is the
stabilizer of the set of vertices $\{v_{d+1},\dots, v_{i}\}$ under $G$,
which we denote by $G_i$.

So, in order to compute the orbits of facets of $p_i(\MP')$ under $G_i$
we need to compute the orbits of facets and ridges of $p_{i-1}(\MP')$,
first under $G_{i-1}$ and then under $G_{i}$ using the double coset
method.

It is well known that the
Fourier--Motzkin method depends in a critical way on the
ordering of generators $v_1,\dots,v_n$ of $\MP$ (see \cite{abs-1997} and \cite{b99}). 
The Cascade Algorithm introduces another such dependency as the size of the
set-stabilizers vary enormously according to the chosen ordering.
Among the many possible orders of the generators,
it is not clear which ordering is best.

In practice Fourier--Motzkin elimination suffers from generating many
redundant inequalities.  This problem can be eliminated by using the
related \emph{Double Description Method}~\cite{fp96}.  For a
description of how to use the double description method to perform the
Fourier--Motzkin projection steps without introducing redundant
inequalities, see~\cite{bremner-thesis}, Section~4.1.

\section{A pivoting method}

\label{sec:sym-pivoting}

Since for certain classes of input, the most successful methods for
the polyhedral representation transformation problem (without taking
symmetries into account) are based on the pivot operation of the
simplex method, it is natural to consider whether pivoting techniques
can be adapted to the symmetric setting.

\subsection{The basis graph up to symmetry}

For the purposes of this discussion, by \emph{$k$-basis}, we mean a
set of $k$ extreme rays generating the linear span of a $k$-face of a
pointed polyhedral cone.  Two $k$-bases are \emph{adjacent} if they
share $k-1$ extreme rays.  When $k$ is not specified, we refer to
$(d-1)$-bases of a $d$-dimensional polyhedral cone.  Note that the
bases of a polyhedral cone depend not only on the combinatorial
information contained in the face lattice, but also on the linear
dependencies among extreme rays incident to a given face.  When
discussing $(d-1)$-dimensional polytopes or their $d$-dimensional
homogenizations, we use basis to refer interchangeably to a
$(d-1)$-basis of the homogeneous cone or a $(d-2)$-basis of the
polytope (i.e., $d-1$ vertices spanning a facet of the polytope).

The \emph{basis graph} has as nodes the bases, and as edges the pairs
of adjacent bases.  In the so-called \emph{non-degenerate} case each facet
of a $d$-dimensional cone is a $(d-1)$-dimensional simplicial cone.
By a perturbation argument, for any cone there exists a triangulation
of the boundary which is combinatorially equivalent to the boundary of
a polyhedral cone with simplicial facets; it follows by Balinski's
Theorem that the basis graph is at least $(d-1)$-connected.

One of the earliest published methods of polyhedral representation
transformation~\cite{charnes53} is based on an exhaustive exploration
of the basis graph.  Using the pivot operation of the simplex method,
adjacent $(d-1)$-bases can be found in time proportional to the number
of input generators and the dimension (for details, see
e.g.~\cite{avis-2000}).  Pivoting methods  based on Reverse
Search~\cite{af92} have memory usage independent of the
output size. In the present work we consider only more direct methods based
on generation (and storage) of the basis graph up to symmetry.

In the typical case, generating the entire basis graph is impractical,
even if it is not stored. This is because of the enormous number of
bases that correspond to each facet in the degenerate (i.e.\
non-simplicial facet) case.
In the non-symmetric setting, perturbation has been widely used to
reduce the size of the basis graph under consideration, 
in particular via the lexicographic (symbolic) perturbation 
discussed further in Section~\ref{sec:symbolic}.  We consider one way in which
perturbation can be applied in the symmetric case in
Section~\ref{sec:orb-perturb}.  
There is some tension between
the notions of symmetry and perturbation though: 
whereas perturbation allows
us to reduce the size of the basis graph,
usually some of the given symmetries are lost. So there is 
a tradeoff in which ideally the quotient graph of the
obtained new basis graph with respect to its remaining symmetries
is as small as possible.

\subsubsection{Symmetry of the basis graph}
The performance of pivoting based methods under symmetry is determined
not by the total size of the basis graph, but by the number of orbits
of bases.  This number is not determined by the number of orbits of
extreme rays and the number of orbits of facets.

We define the \emph{basis automorphism group} to be the subgroup of
the combinatorial automorphism group that acts on the basis graph.  To
see that this is a non-trivial restriction, consider the cone
generated by the following (row) vectors in $\R^5$
\begin{equation}
  \label{eq:oct-pyr}
  \begin{matrix}
    1&0&1/2&1 & 1\\
   -1&0&1/2&1 & 1\\
    0&1&0&1 & 1\\
    0&-1&0&1 & 1\\
    0&0&1&1 & 1\\
    0&0&-1&1 & 1\\
    0&0&0&-1 & 1
  \end{matrix}
\end{equation}
Combinatorially, it is (the homogenization of) a pyramid over an
octahedron, and its combinatorial automorphism group is the $48$
element  octahedral group.  The basis automorphism
group has only $16$ elements; in particular the orbit of the first $d$
generators (as a basis) is of size $1$.  

\medskip

\begin{example}  \label{ex:crosspolytope}
  The regular $d$-cross polytope $C_d$ has 
  one orbit of $d$-bases and one orbit of $(d-1)$-bases
  with respect to its automorphism group.
\end{example}

\begin{proof}
  Since the $d$-cross polytope is simplicial, 
  (i.e., its facets are all simplices)
  its $(d-1)$-bases are exactly its facets. Observe that
  at most one pair of vertices in a $d$-basis must consist of an
  opposite pair $\pm e_j$, since any two pairs are co-planar.  It thus
  follows from the pigeonhole principle that a $d$-basis consists of a
  facet, along with one vertex not on that facet.  All such simplices 
  are equivalent under the automorphism group of the cross polytope.
\end{proof}

    For a slightly more involved example,
    let $\MP = \conv(v_1, \dots, v_m) \subset \R^d$
    and $\MQ = \conv(w_1, \dots, w_n) \subset \R^e$. 
    The \emph{wreath product}  of $\MP$ with $\MQ$ is
    defined as 
    \begin{equation*}
      \MP \wr \MQ := \conv \{ ( \underbrace{0,\dots ,0}_{d(k-1)},
      v_i,\underbrace{0,\dots ,0}_{d(n-k)},w_k)  : 
            1 \leq i \leq m , 1 \leq k \leq n \}
      \subset \R^{nd+e}. 
    \end{equation*}

    So loosely speaking, the wreath product of $\MP$ and $\MQ$ is
    obtained by attaching to the vertices of $\MQ$ pairwise  
    orthogonal subspaces, each containing a copy of $\MP$.

    \begin{example}
      The wreath product $\MP\wr \MQ$ of a regular 
      $d$-cross polytope $\MP$ with a regular $e$-cross polytope $\MQ$ has
      dimension $D=2de+e$, $4de$ vertices, $2^{(d+1)e}$ facets
      and one orbit of vertices, facets, and $(D-1)$-bases
      with respect to its linear automorphism group. 
      \label{ex:wreath}
    \end{example}

    \begin{proof}
      The proof is based on the results of Section~2.2 in~\cite{jl05}.
      The number of vertices of the cross polytopes $\MP$ and $\MQ$ are
      $m=2d$ and $n=2e$. The assertion on the dimension and 
      the number of vertices of $\MP\wr \MQ$ follows immediately from the 
      definition of wreath products.  
      The number of facets and the number of bases orbits can
      be derived from Proposition~2.2 in~\cite{jl05}.
      By it, the facets of $\MP\wr \MQ$ are in one-to-one correspondence
      with choices $(F;F_1,\dots, F_e)$ of a facet $F$ of $\MQ$ and
      facets $F_i$ of $\MP$; which gives the count of facets.
      Assume w.l.o.g. that $F$ contains the vertices                
      $v_{1},\dots, v_{e}$ of $\MQ$. Then the rows of the following
      matrix describe the vertices contained in a facet of $\MP\wr \MQ$:
      \begin{equation*}
        \begin{bmatrix}
          \MP & &      &    &     &     &        &    & v_{1}\\
            &\MP&      &    &     &     &        &    & v_{2}\\
            & &\ddots&    &     &     &        &    &      \\
            & &      & \MP  &     &     &        &    & v_{e}\\
            & &      &    & F_1 &     &        &    & v_{e+1}\\
            & &      &    &     & F_2 &        &    & v_{e+2}\\
            & &      &    &     &     & \ddots &    &\\
            & &      &    &     &     &        & F_e& v_{2e}
        \end{bmatrix}
      \end{equation*}
      The entries $\MP$ and $F_i$ stand for all choices of vertices
      from $\MP$, respectively $F_i$. So a facet of $\MP\wr \MQ$ contains
      $3de$ vertices.

      The linear automorphism group of $\MP \wr \MQ$ contains
      the {\em semidirect product} $G=G_P^n \rtimes G_Q$
      (where $G_P$ and $G_Q$ denote the automorphism
      groups of $\MP$ and $\MQ$ and $G_P^n$ denotes an $n$-fold 
      direct product). In particular, $G_Q$ is
      isomorphic %
      to a subgroup of $G$, which permutes not only $v_1,\dots,v_{2e}$ 
      in the last $e$ coordinates, but also the $n$ copies
      of $\MP$ attached to them. 
      In contrast to this action interchanging the
      $n$ \ $d$-dimensional subspaces, the elements of the %
      subgroup~$G_P^n$ of~$G$ act only within these subspaces.
      By this, all of the vertices and facets lie in one orbit 
      under the action of~$G$.
      The assertion on the bases follows as well,
      after noting that a base within a facet (as described above)
      is determined by $e$ choices of $d$-bases 
      (in the $e$ copies of $\MP$). All of them are equivalent
      with respect to $G$ by Example \ref{ex:crosspolytope}.
    \end{proof}

\subsection{A prototype implementation}

We have implemented a prototype called \emph{symbal} (``Symmetry,
Bases, and Lexicography'') for facet enumeration up to symmetries via
pivoting.  The code is mainly in GAP~\cite{GAP}, with an external
server based on \texttt{lrs} doing the pivoting. The main algorithm is
a depth first search of the pivot graph; several options are provided
for symmetry testing and perturbation (discussed further below).  
All symmetry groups are represented internally as permutation groups,
and  the 
backtracking algorithm of GAP's \textsc{RepresentativeAction} is 
used to test bases and facets for equivalence.
Our
main goal with this implementation is not to be able to attack large
problems, but to provide a test-bed for better understanding of the
main issues involved with pivoting under symmetry.  
In general we thus are 
not too concerned with a constant multiplicative slowdown in the
runtime of our code versus more polished software; we are more
interested in what kinds of problems can be solved in a reasonable
amount of time.

\begin{figure}
  \centering
  \includegraphics[height=2.5in]{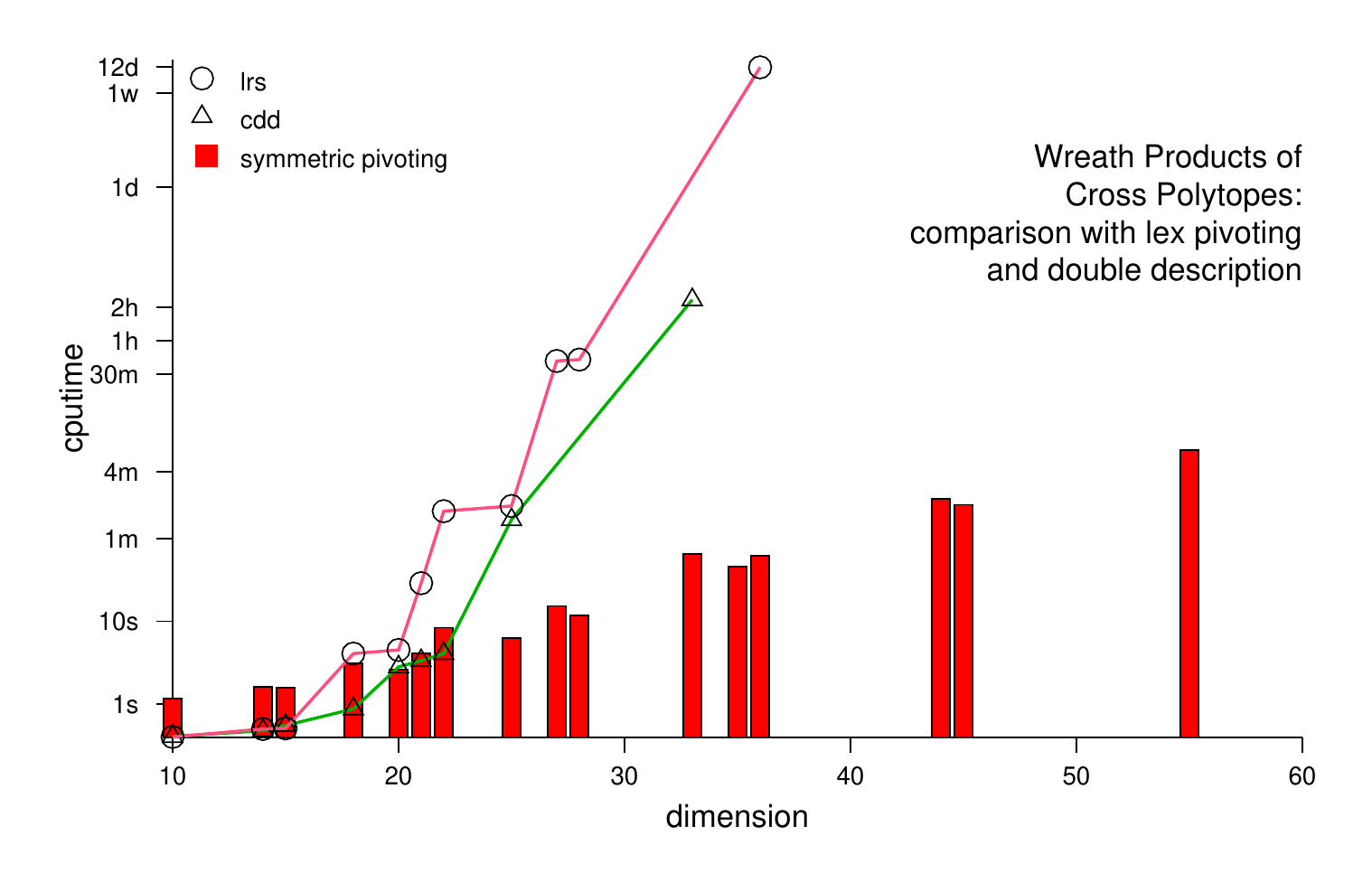}
  \caption{Experimental results for wreath products of cross polytopes}
  \label{fig:wreath-compare}
\end{figure}

At least for certain special cases, our prototype is able to solve
quite large problems.  Figure~\ref{fig:wreath-compare} compares the
run time of \texttt{symbal} with \texttt{cdd} and \texttt{lrs} on the
wreath products of cross polytopes discussed in
Example~\ref{ex:wreath}.  Of course both of \texttt{cdd} and
\texttt{lrs} compute all of the facets rather than orbits. 
The times reported here are on a 3GHz
Pentium IV with 1G of memory.

We call a cone (or polytope) \emph{basis-simplicial} if it has as few
basis orbits as facet orbits. Evidently every simplicial polytope is
basis-simplicial.  For perhaps less contrived examples of
basis-simplicial polytopes than Example~\ref{ex:wreath} we mention the
Dirchlet-Voronoi-cells (DV-cells) of the root lattices ${\mathsf
  D}_3$, ${\mathsf D}_4$, and ${\mathsf E}_8$ (convex hull of their
shortest non-zero vectors) and some related {\em universally optimal}
spherical polytopes (see \cite{universal-optima}).  Their boundaries
consist of regular cross polytopes and regular simplices only and they
thus have at most two orbits of basis with respect to their symmetry
group (cf.\ Example~\ref{ex:crosspolytope}).

As with the non-symmetric case, for cones with a large number of
(orbits of) bases, pivoting is not a good approach, at least not
without some kind of perturbation.  One benchmark of how 
{\em orbitwise degenerate} a polytope  is, 
is whether the number of orbits of bases is
larger  than the number obtained by triangulating the boundary. 
In this case,
one is almost certainly better off applying either \texttt{cdd} or
\texttt{lrs}.  A familiar example 
where this level of degeneracy occurs is in the $d$-dimensional cubes
(see Table~\ref{tab:cube-orbits}).
    \begin{table}
  \centering
    \begin{center}
      \begin{tabular}{ccccc}
        Dimension & \texttt{lrs} Triangulation  &  Basis Orbits \\\hline\hline
        4 & 48 & 4 \\ 
        5 & 240 & 17 \\
        6 & 1440 & 237 \\
        7 & 10080 & 9892 \\
        8 & 80640 & $>209000$ 
    \end{tabular}
    \end{center}
  \caption{Basis orbits of cubes}
  \label{tab:cube-orbits}
\end{table}

\subsection{Orbitwise Perturbation}
\label{sec:orb-perturb}

The basis automorphism group could in principle be larger than the
linear automorphism group (as a simple example, consider a simplicial
polytope combinatorially equivalent to the regular $d$-cross polytope,
but with trivial linear symmetry group).  Nonetheless, for the reasons
articulated in Section~\ref{sec:symmetries}, and because linearity
simplifies the discussion here, we consider here the case where we are
given a subgroup of the restricted isomorphism group of a polyhedral
cone.

In this section we consider modifications of the standard 
\emph{lexicographic perturbation} schemes that
preserve some, but not necessarily all of the symmetry of the input.
We first  characterize one kind of transformation that preserves 
a prescribed symmetry group.
\begin{proposition}
  \label{prop:fixedpoint}
  Let $V$ be a vector family with restricted automorphism group $G$.
  Let $H \leq G$.
  Let $V_1, V_2, \dots , V_k$ be the orbits of $V$ under $H$. 
  Let
  $u_1, \dots , u_k$ be a family of fixed points for $H$. Let 
  \begin{equation*}
    V'= \bigcup_j \{ v_i +u_j \mid v_i \in V_j\}
  \end{equation*}
  Let $H'$ be the restricted automorphism group of $V'$.
  Then $H \leq H'$.
\end{proposition}

\begin{proof}
  Let $T\in H$ be a restricted automorphism for $V$ and $\rho$ the
  corresponding permutation such that $Tv_i=v_{\rho(i)}$.  Suppose
  $v_i \in V_j$.  Since there exists an automorphism mapping $v_i$ to
  $v_{\rho(i)}$, by definition $v_{\rho(i)} \in V_j$.
  It follows by linearity of $T$ that 
  $$
      T(v'_i) = T(v_i+u_j) = Tv_i+Tu_j = v_{\rho(i)}+u_j = v'_{\rho(i)}.
  $$
\end{proof}

Of course, in addition to preserving part of the symmetry group, we
need to solve the original representation conversion problem.  Let
$\nu(\cdot)$ denote the orbitwise perturbation map.  The key property
we need is that for each facet $F$ of the original cone, there is some
basis $B$ such that $\nu(B)$ is a basis of the perturbed cone.  The
standard way of ensuring this, and the one we adopt here, is to insist
that the facets of the perturbed cone induce a subdivision of the
facets of the original cone.  Several notions of perturbation exist in
the literature, including lexicographic~\cite{DantzigOrdenWolf55},
numeric~\cite{MegChand89}, symbolic~\cite{yap90}, and
geometric~\cite{seidel98}.  We will discuss a method that can be seen
as a \emph{linear perturbation} in the language of
Seidel~\cite{seidel98} or a modified version of the lexicographic
perturbation first proposed by Dantzig, Orden and
Wolf~\cite{DantzigOrdenWolf55}.

Here we take a somewhat weaker definition of perturbation than is
typical, since we are not concerned necessarily with obtaining a
cone with simplicial facets as result.  Let $V$ be vector family. We say
$W\subseteq V$ is \emph{extreme} (for $V$) if it is contained in some
facet of $\cone(V)$.  
\begin{definition}
\label{def:validperturb}
We say that $\widetilde{V}$ is a \emph{valid
  perturbation} of $V$ if there is a bijection $\nu(\cdot)$ between
$V$ and $\widetilde{V}$ such that for any $W\subseteq V$,
\begin{enumerate}
\item  If $\nu(W)$ is linearly dependent then $W$ is.
\item  If $\nu(W)$ is extreme for $\widetilde{V}$ then $W$ is extreme for $V$.
\end{enumerate}
\end{definition}

From Definition~\ref{def:validperturb}, we have immediately the
following refinement property. 

\begin{proposition}
  \label{prop:valid1}
For a valid perturbation $\widetilde{V}$ of a vector family $V$ we have
  \begin{enumerate}
  \item
    The boundary complex of $\cone(\widetilde{V})$ 
    is a subdivision of the boundary complex of $\cone(V)$.
  \item
    If $\widetilde{X}$ is a valid perturbation of
    $\widetilde{V}$, then $\widetilde{X}$ is a valid perturbation of $V$.
  \end{enumerate}
\end{proposition}

The following proposition shows that valid perturbations  result
from sufficiently small changes to a vector family.

\begin{proposition}
  \label{prop:valid2}
  For any vector family $V\subset \R^d$, for any $W \subseteq V$ and any
  vector $u\in \R^d$, 
  $\widetilde{V}(\varepsilon)=\{\, w+ \varepsilon u : w \in W\} \cup V\setminus W$ 
  is a valid perturbation of $V$ for all $\varepsilon$ with 
  $|\varepsilon|$ sufficiently small.
\end{proposition}
\begin{proof}  
  Consider the vector family $\widetilde{V}(\varepsilon)$ changing as
  $\epsilon$ varies. From the point of view of linear dependence,
  we need only concern ourselves with a $d$-set of vectors $v_1, \dots, v_d$, 
  which are independent in $V$, but whose perturbation becomes
  dependent for $\epsilon>0$. 
  The following continuity argument shows that for small enough $\varepsilon$
  no new linear dependencies may occur.
  Let $T$ be the matrix whose rows are
  $v_1, \dots , v_d$.  By renumbering we may assume the first $k$ rows of
  $T$ are in $W$.
  Let $U$ be the matrix whose first
  $k$ rows are $u$ and the remaining rows are the zero vector.  
  Then
  \begin{equation*}
    g(\varepsilon)=\det (T+\varepsilon U)
  \end{equation*}
  is a polynomial in $\varepsilon$, which is non-zero in a
  neighborhood of $\varepsilon=0$.
  Taking the intersection of all such
  intervals over all full rank $d$-sets of vectors yields 
  (1) and (2) of Definition~\ref{def:validperturb}
\end{proof}

Let $V$ be a vector family with restricted automorphism group $G$.
Let $H$ be a subgroup of $G$.  Let $u_H$ be a fixed point for $H$.
We say that $\widetilde{V}$ is obtained by \emph{pushing} $W\subseteq V$
(respectively
\emph{pulling} $W\subseteq V$) if
\begin{equation*}
  \widetilde{V} =  \{\, w +\sigma \varepsilon u_H \mid w \in W\,\}
  \cup V\setminus W
\end{equation*}
where $\sigma=1$ (respectively $\sigma=-1$) and $\varepsilon>0$ is
sufficiently small so that $\widetilde{V}$ is a valid perturbation of
$V$.  

Let $V_1, V_2, \dots V_k$ be the orbits of $V$ with respect to $H$. In general the
combinatorial structure of the resulting boundary complex depends not
only whether each $V_j$ is pulled or pushed, but on the order these
operations are carried out.
We say that $\widetilde{V}$ is an \emph{orbitwise lexicographic perturbation} of $V$ with respect to $H$
if it is obtained by  pulling or pushing each orbit defined by $H$ in
some \emph{perturbation order} $\pi$.
  Considering the special case where $V$ contains
the homogenization of the vertices of a $d$-polytope 
and $u_H$ is the homogenization of the origin, 
pulling (respectively pushing) corresponds to scaling
an orbit outward (respectively inward) with respect to the origin;
this terminology is consistent with that of Lee~\cite{lee91}.

Combining Propositions~\ref{prop:fixedpoint}, \ref{prop:valid1},
and~\ref{prop:valid2}, we have:
\begin{proposition}
  If $\widetilde{V}$ is an orbitwise lexicographic perturbation of $V$ 
  with respect to a subgroup $H$ of $V$'s restricted automorphism group, 
  then $H$ acts on the basis graph of $\cone(\widetilde{V})$.
\end{proposition}

In the case that there is only one orbit under the basis
automorphism group, perturbing all input vectors by the same vector
will not decrease the degeneracy of the problem.  Our general strategy
will thus be to choose a subgroup $H$ that has multiple orbits of
extreme rays, and at the same time choose different perturbations $u_i$ for each
orbit of $V$ with respect to $H$.

It is known~\cite{lee91} that lexicographic perturbation of all of the vertices
induces a triangulation (also called lexicographic) of the
boundary of the polytope.  In our case because we perturb (i.e.\ push
or pull) all vertices in an orbit by the same amount, we cannot in
general guarantee a simplicial result (i.e.\ an induced triangulation).
Thus we continue to explore the complete basis graph of the perturbed
cone (polytope), up to equivalence classes of bases.  We nonetheless
hope for a significant reduction in the size of the basis graph by
(effectively) breaking up very degenerate facets of the cone.  The
tradeoff, examined further in the next section, is that 
although we may lose some symmetry of the polyhedron, the 
quotient graph of the new basis graph (with respect to the remaining symmetries)
may become smaller. %
Moreover, the degree of the vertices in the basis graph
may have decreased, which speeds up the computation of the quotient.

\subsubsection{Choosing a subgroup to preserve}

For any polyhedral cone $\MP$ with symmetry group $G$, in order to
effectively use orbitwise perturbation, one needs to find a subgroup
$H\leq G$ such that the above mentioned possible computational gain is
as large as possible. In order to develop some heuristics for how one
might find such a subgroup, we have systematically studied the
subgroups of the restricted automorphism group of the DV-cell
of the $\SE_7$ root lattice (hereafter we use $\conv \SE_7$ to denote
this DV-cell).  Although we do not claim the results from a
single example are in any way conclusive, they do at least suggest
some ideas for further study.

Let $G$ denote the restricted automorphism group of $\conv \SE_7$.
This is a group of order $2903040$ that has one orbit on the $126$
vertices of $\SE_7$. The polytope $\conv \SE_7$ has $632$ facets in
$2$ orbits. It has $161$ basis orbits, compared to $20520$ bases in
the triangulation produced by \texttt{lrs}.

The experiments on this section were carried out an Acenet
(\url{http://www.ace-net.ca}) cluster with SunFire x4100 nodes (two
2.6 GHz dual-core Opteron 285 SE processors and 4 GB RAM per core).

\begin{figure}
  \centering
  {\includegraphics[width=\textwidth]{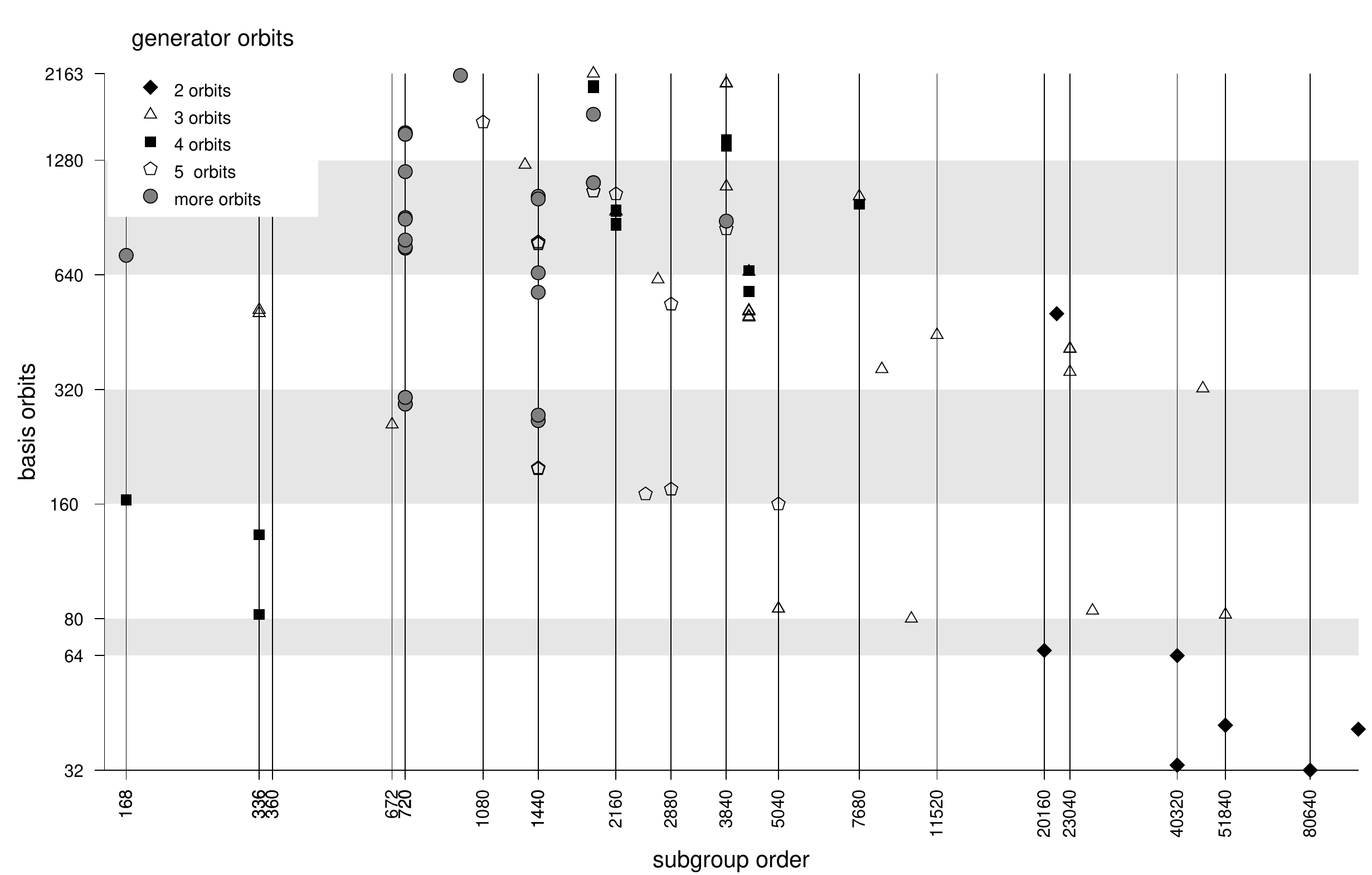}}
  \caption{Number of basis orbits for various subgroups on  $\conv {\mathsf E}_7$}
  \label{fig:E7-bases}
\end{figure}

\begin{figure}
  \centering
  {\includegraphics[width=\textwidth]{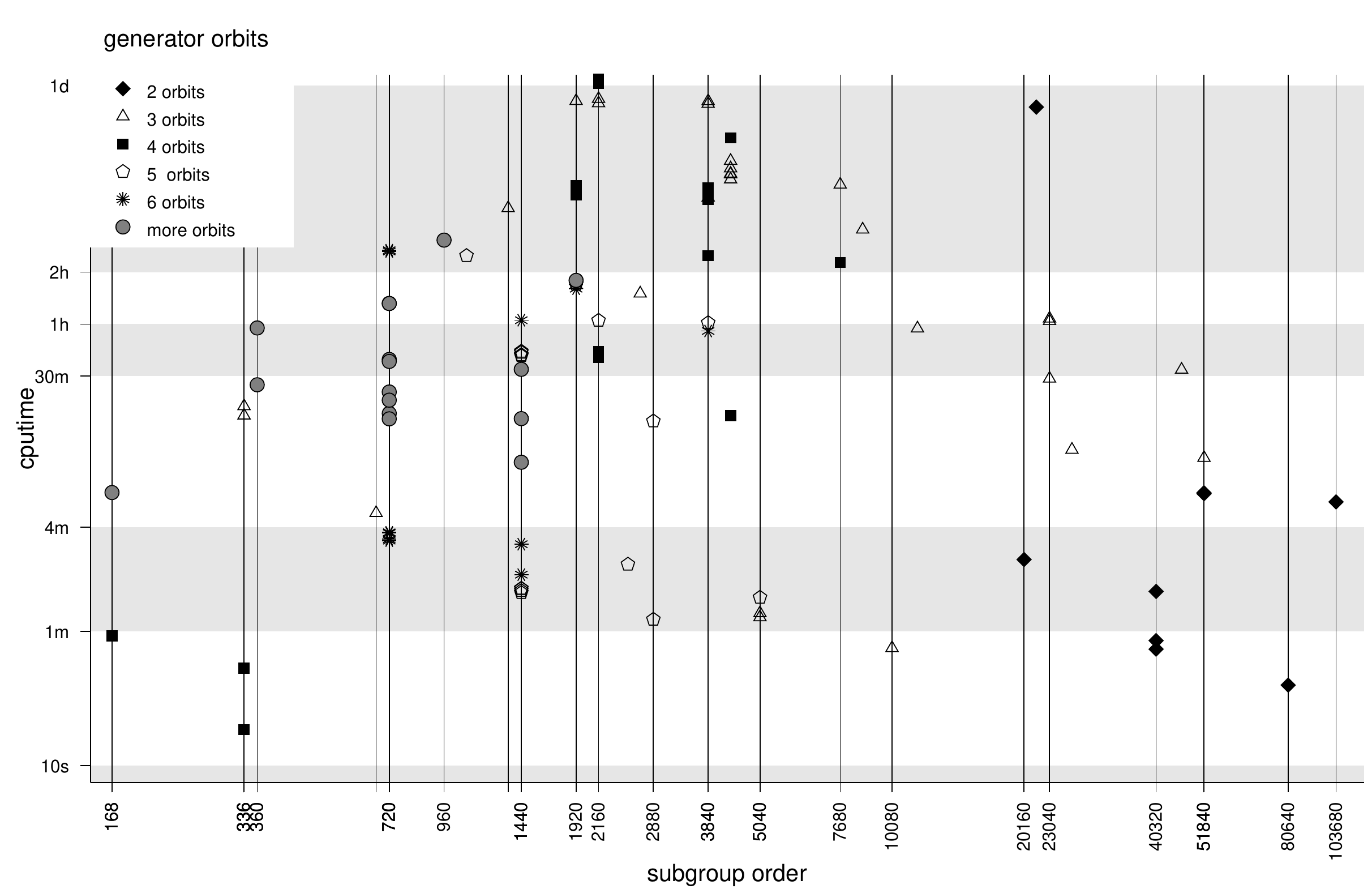}}
  \caption{CPU time to enumerate facets for various
    subgroups on  $\conv {\mathsf E}_7$}
  \label{fig:E7-time}
\end{figure}

A sample of conjugacy classes of the subgroup lattice was generated
using the GAP function {\tt LatticeByCyclicExtension}, with the
restriction that groups of size less than or equal to 100 were
discarded. A representative was chosen from each conjugacy class,
yielding $102$ subgroups of $G$.  In general combinatorially distinct
perturbations could result from choosing a different perturbation
order, and by varying the direction (push or pull) that each orbit is
perturbed.  In this experiment we restricted ourselves to orbitwise
pulling, and ordered according to the smallest index in an orbit.
Orbitwise pulling perturbations were generated for each of these
subgroups, and these perturbations were used as input to our prototype
\texttt{symbal}. A total of $99$ of these computations completed in
the time allocated.  The fastest time was about 16 seconds, and
longest just over $26$ hours. The unperturbed version completed in
about half an hour. Perturbation thus yielded results ranging from
more than a $100$-fold speedup to more than $50$-fold slow down.

Figure~\ref{fig:E7-bases} presents the experimental results in terms
of the number of basis orbits computed. The groups are classified
according to their order ($x$-axis) and number of orbits of vertices (
shown by symbol).  Recall that for any actual perturbation, we need at
least two generator orbits.  Of the subgroups with less than $64$ orbits
of (perturbed) bases all had two generator (i.e.\ vertex) orbits.  The
best result is $32$ basis orbits, achieved by a subgroup of order
$80640$. Within the class of subgroups having two generator orbits,
the best results were achieved by groups of large order. Note that a
large subgroup inducing two (or any other fixed small number) of
generator orbits can be found by choosing random sets of elements as
generators.

Comparing Figure~\ref{fig:E7-bases} with Figure~\ref{fig:E7-time}, we
see runtime is not completely a function of the number of basis
orbits, but the size of the group also plays a role.  The actual best
time is achieved by a group of order $336$, with $4$ generator orbits and
$82$ basis orbits.  After a certain point a small subgroup is no
longer advantageous, since the number of basis orbits are too large. In
the limiting case, with a trivial group, our pivoting scheme produces
a lexicographic triangulation like \texttt{lrs}; albeit less efficiently
in our prototype implementation.  
A relatively small subgroup with a fixed
number of input orbits can also be found by the same random sampling
strategy. 

\begin{figure}
  \centering
  \begin{tabular}{c@{\hspace*{0.2in}}c@{\hspace*{0.2in}}c}
  \includegraphics[width=0.3\textwidth]{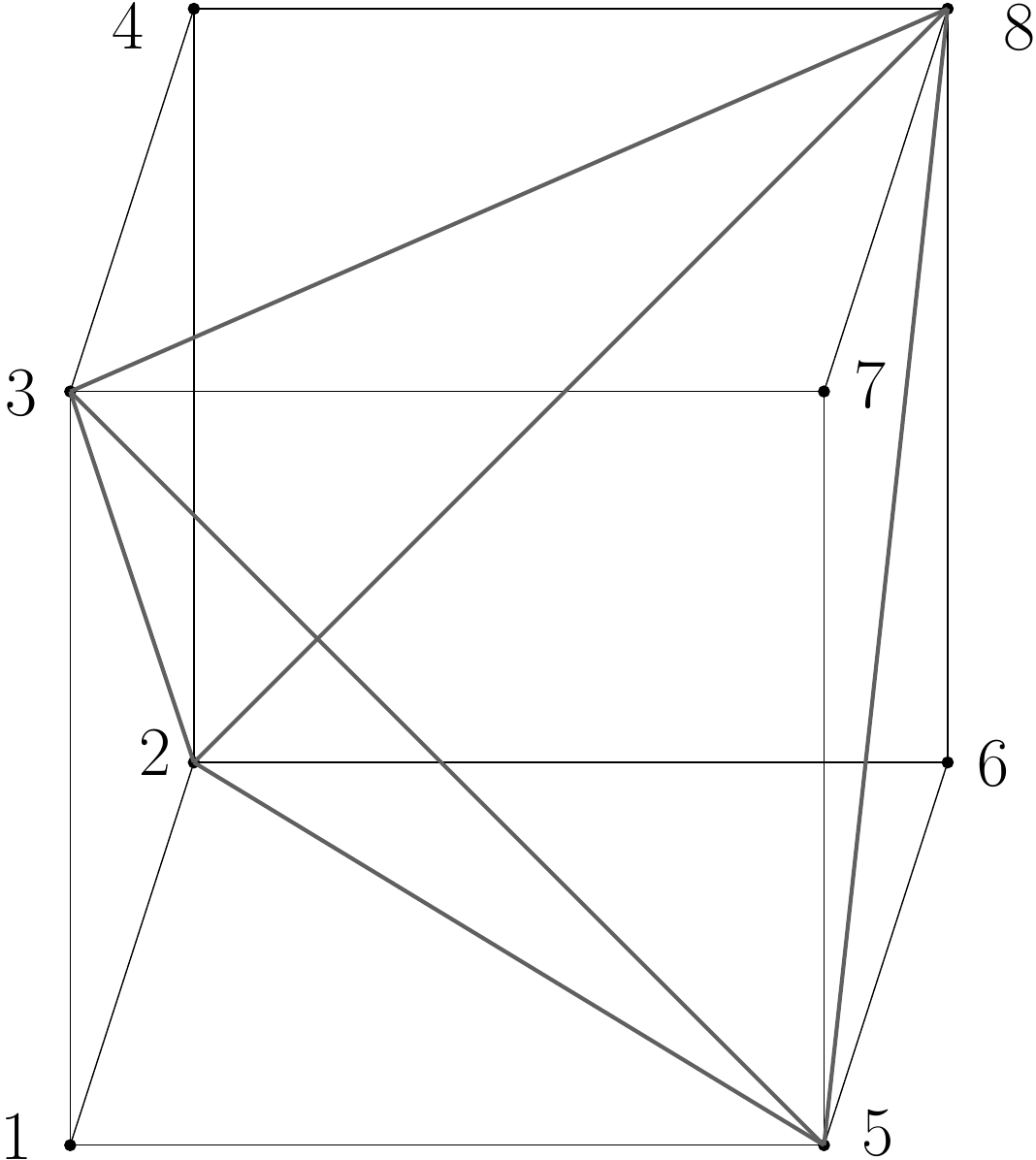} &
  \includegraphics[width=0.3\textwidth]{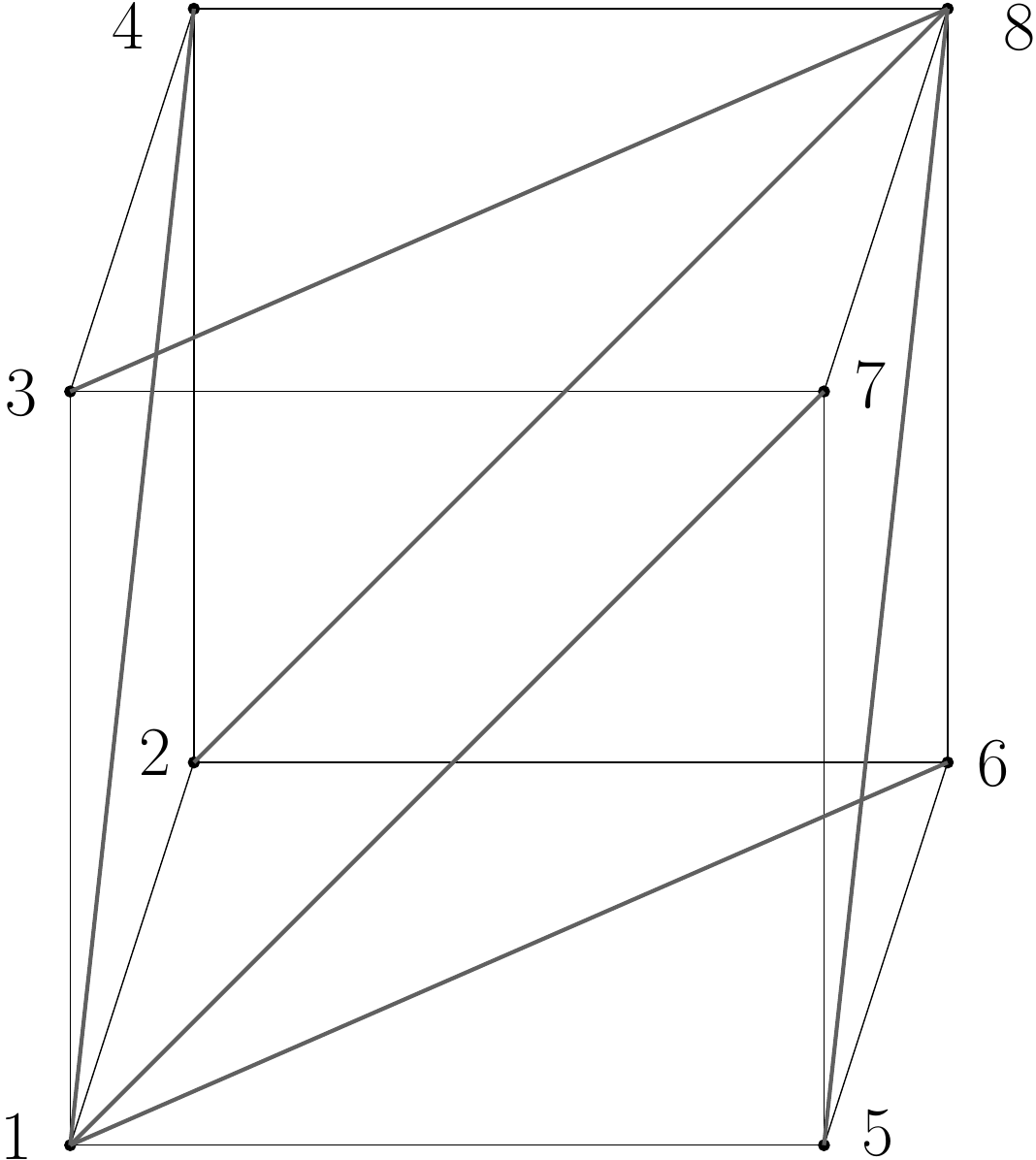}&
  {\includegraphics[width=0.3\textwidth]{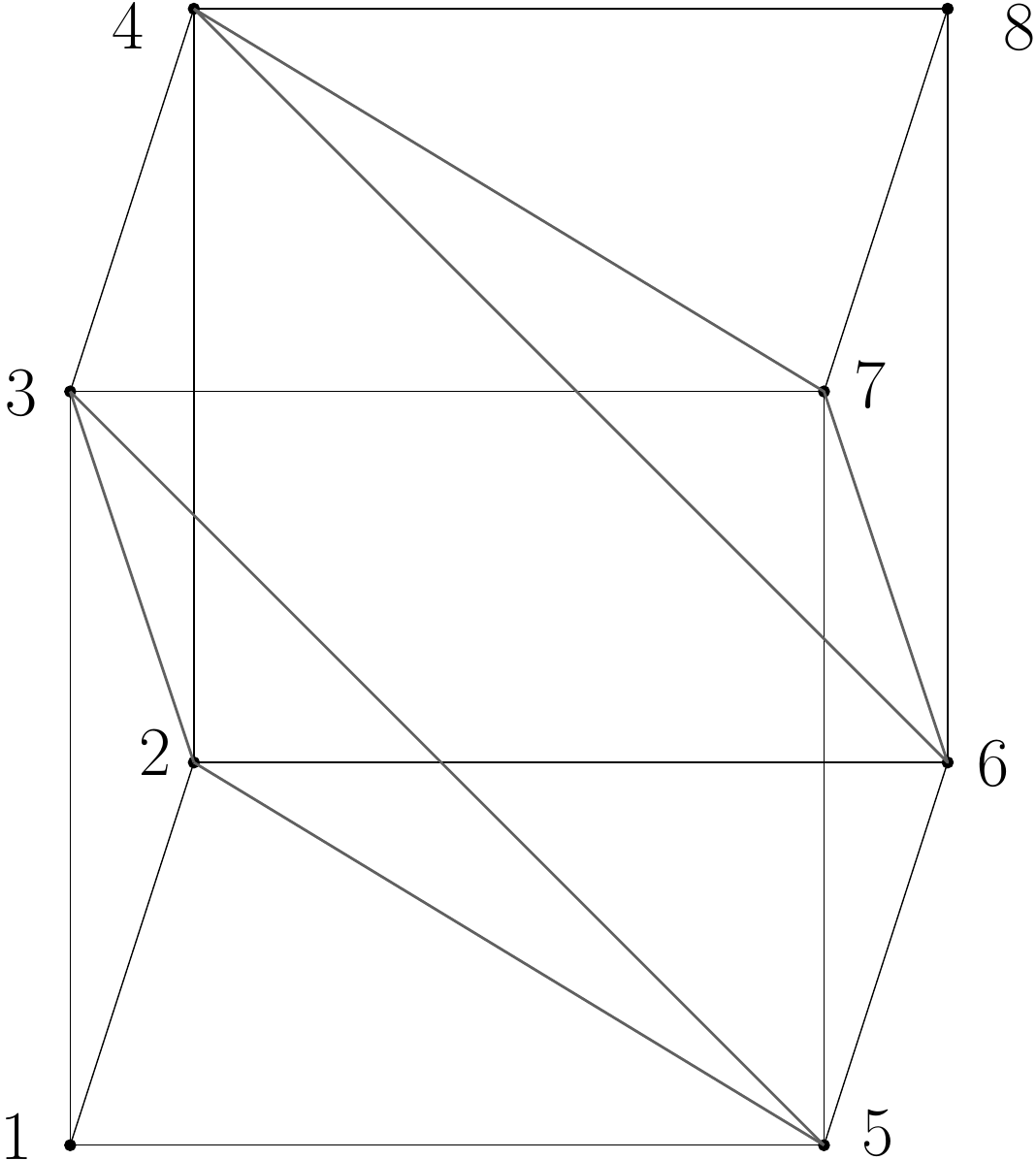}}\\
  $\pull(2,3,5,8)$& $\pull(1,8)$ & $\push(1,8)$
  \end{tabular}
  \caption{Symmetric triangulations of the boundary of the $3$-cube; 
  in each case one orbit push or pull suffices to triangulate.}
  \label{fig:cube-tri}
\end{figure}

Figure~\ref{fig:cube-tri} illustrates our two heuristics applied to
the three dimensional cube. On the left we have the triangulation
induced by a subgroup of order $24$ having two generator orbits, 
In the center, the triangulation is induced by the stabilizer 
of a pair of opposite vertices.  On the right, we take the same
stabilizer, but push the first orbit instead of pulling.

The triangulation induced by pulling opposite vertices of a $3$-cube
turns out to be a special case of a triangulation of the $d$-cube that
has only one basis orbit in any dimension.  Let $I^d=[-1,1]^d$ denote
the centrally symmetric $d$-cube.  Define~$\one=\sum_{i=1}^{d} e_i$
For each permutation $\rho\in \Sym(d)$, let
$\Delta_\rho$ denote the simplex with vertices:
\begin{equation*}
  -\one,\; -\one+2(e_{\rho(1)}),\; -\one+2( e_{\rho(1)}+e_{\rho(2)}),\dots,\;\one
  \end{equation*}
  Let $\Delta$ denote the union of all $\Delta_\rho$ and let
  $\bar{\Delta}$ denote the set of $(d-1)$-simplices formed by
  intersecting the simplices of $\Delta$ with the boundary of $I^d$.
  It is known~\cite{drs-2007} that $\Delta$ forms a triangulation 
  of $I^d$; consequently $\bar{\Delta}$ forms a triangulation of the
  boundary, since every triangulation of a polytope induces a
  triangulation of its boundary.  We call $\Delta$ (respectively
  $\bar{\Delta}$) the \emph{linear ordering triangulation} of $I^d$
  (respectively of the boundary of $I^d$.)  Since there is a bijection
  between permutations in $\Sym(d)$  and simplices in $\Delta$, 
  $\Sym(d)$ acts transitively on
  $\Delta$ by permuting coordinates.

\begin{example}
  Let $V$ denote the vertices of $I^d$, i.e.\ $V=\set{\pm 1}^d$.
  Let $H$ denote the stabilizer of the  automorphism group of $I^d$ 
  on $\set{ -\one, \one}$. 
  Define $\omega(v)=\min(\one^Tv,-\one^Tv)$.
  Let $\nu(V)$ denote the $H$-orbitwise pulling of $V$ in the order induced 
  by $\omega$. Then the following holds:
  \begin{enumerate}
  \item[(a)] $H$ acts transitively on $\bar{\Delta}$.
  \item[(b)] $\partial \conv\nu(V)$ is combinatorially equivalent to
    $\bar{\Delta}$.
  \end{enumerate}
\end{example}
\begin{proof}
  Let us first remark that the order of orbits induced by $\omega$ is
  well defined.  The group $H\leq \Sym(2^d)$ is generated by a set of
  generators permuting coordinates, along with a \emph{switching
    permutation} $\sigma$ that maps $v$ to $-v$.
  It follows that the $H$-orbits of $V$
  are the equivalence classes of $V$ with respect to $\omega$.
  
  To see (a), consider the $d$-simplex (corresponding to the identity
  permutation) $\Delta_{\Id}$ with vertices $\set{ -\one, -\one+2e_1, -\one+2(e_1+e_2),\dots, \one}$.
 Let $\sigma'$ denote the permutation in $H$
  that first applies $\sigma$ (i.e.\ switching) followed by reversing
  the order of coordinates. The permutation $\sigma'$ is an
  automorphism of $\Delta_{\Id}$ which carries the $(d-1)$-simplex
  $\set{-\one, -\one+2e_1,\dots,e-2e_d)}$ to  $\set{ \one,  \one-2e_d,  \dots, -\one + 2e_1}$.  
These two $(d-1)$
  simplices are precisely the contribution of $\Delta_{\Id}$ to the
  linear ordering triangulation $\bar{\Delta}$.  Thus any simplex of
  $\bar{\Delta}$ can be mapped to any other by an action of $\Sym(d)$
  on $\Delta$ (i.e.\ permuting coordinates), followed by possibly
  applying $\sigma'$.

  We now consider (b).  We argue that $\nu(V)$ induces a linear
  ordering triangulation of each $k$-face of $I^d$, $0\leq k < d$.
  Each $(d-2)$-face will receive the same (linear ordering)
  triangulation from the two facets that contain it, hence the
  triangulations of the facets form a triangulation of the boundary.

  For $0 \leq k \leq 2 $, there is nothing to prove.  Let $F$ be a
  $k$-face of $I^d$, $2<k<d$.  Let $v^+$ (resp.\ $v^-$) be the vertex of
  $F$ with the most positive (resp.\ negative) coordinates.  Recall
  that we will first pull the $H$-orbit with smallest $\omega$ value.

  If the functional $\omega$ is minimized uniquely at $v^*\in
  \set{v^-, v^+}$ then the perturbation corresponds locally to a
  standard~\cite{lee91} pulling and the corresponding subdivision is
  into pyramids $F_1 ,\dots, F_j$ with apex $v^*$ and bases
  corresponding to all of the $k-1$ faces of $F$ that do not contain
  $v^*$.

  Otherwise $\omega$ is minimized at both $v^-$ and $v^+$.  The
  perturbation thus takes $v^+$ to $\rho v^+$ and $v^-$ to some $\rho
  v^-$, $\rho > 1$. This turns out to induce a subdivision of $F$ into
  polytopes $F^2_j$ with vertices $V_j=\set{ \rho v^-, \rho v^+ } \cup
  R_j$ where $R_j$ is the vertex set of a $(k-2)$-face of $F$
  containing neither $v^-$ nor $v^+$.  The polytope $F^2_j$ is a
  \emph{ $2$-fold pyramid}, since $v^- \notin \aff R_j$ and $v^+
  \notin \aff (\set{v^-} \cup R_j)$.  It follows that $\dim F^2_j =
  \dim F=k$.

  That the $F^2_j$ are induced by the perturbation can be seen 
  by exhibiting a supporting hyperplane of $\conv \nu(V)$.
  Without loss of generality, let $F$ be defined by equations
  $x_i=1$, $i=1,\dots, d-k$.  The $(k-2)$-face $\conv(R_j)$ must be defined by further
  equations $x_p=1$, $x_q=-1$, $p,q > d-k$.  Let
  $\mu=1/\rho$. Consider the hyperplane $h_j=\set{ x \mid a^Tx=1 }$,
  where
  \begin{equation*}
    a_i=
    \begin{cases}
      \mu/(d-k) & 1\leq i\leq d-k\\
      (1-\mu)/2 & i=p\\
      -(1-\mu)/2 & i=q\\
      0 & \text{otherwise}\,.
    \end{cases}
  \end{equation*}
  It can be verified that $h_j$ supports $\conv \nu(V)$ and $h_j \cap
  \conv\nu(V) = F^2_j$.

  It remains to see that the $F^2_j$ cover $F$, i.e.\ that there are
  no other cells in the induced subdivision.  Consider an arbitrary
  relative interior point $x$ of $F$.  Let $r$ be the ray from $\rho v^-$
  through $x$. Let $y$ be the first intersection of $r$ with $\partial
  F$ after $x$.  For each $(k-1)$-face of $F$, the double pulling of
  $\set{v^+,v^-}$ acts like a single pulling decomposing the boundary
  of $F$ into pyramids with apex either $\rho v^-$ or $\rho v^+$; hence $y \in
  \partial F^2_j$ for some $j$.  It follows  that $x$ is in $F^2_j$.

  Now suppose for all $j<k$, $\nu(V)$ induces a linear ordering
  triangulation of the $j$-faces of $I^d$.  From the refinement
  property Proposition~\ref{prop:valid1}, we know that $\nu(V)$
  induces a decomposition of the pyramids $F_1 ,\dots, F_k$
  (respectively of the $2$-fold pyramids $F^2_1, \dots,
  F^2_{k(k-1)}$).  In both cases the resulting $k$-simplices
  correspond to the coordinatewise-monotone paths from $v^-$
  to~$v^+$ in $F$.
\end{proof}

\subsubsection{Symbolic implementation}
\label{sec:symbolic}

Suppose we are given $V'\subset\R^{d+1}$ which are homogeneous
coordinates for the generators $V$ of some polyhedron $\MP\subset
\R^d$. We  consider $V$ here as a matrix (with $v_i$ as rows), and let
$u$ denote the column vector of corresponding $(d+1)$st coordinates in
$V'$.  Following the conventions of Section~\ref{sec:basic_notations},
we suppose $u_{j}$ is $1$ if $v_j$ is a vertex, and $0$ if it is an
extreme ray.  We consider the polyhedron $\MP^\diamond$
\begin{equation}
\label{eq:mpdiamond}
\MP^\diamond  = \{\, x\in\R^d : Vx \geq -u \,\}\,.
\end{equation}
The polyhedron $\MP^\diamond$ may be thought of
as $\cone(V')^*\subset (\R^{d+1})^{\ast}$ intersected with the
hyperplane $x_{d+1}=1$. By duality, 
to find the generators of $\MP^\diamond$ is equivalent to 
finding the facets of $\MP$.

Let $G$ be the restricted automorphism group of $V$. Let $G'$ be the
induced group acting on $V'$.  We will assume without loss of
generality that the origin is the centroid of V and thus a fixed point
of $G$. It follows that $e_{d+1}$ is a fixed point of $G'$.  Applying
Proposition~\ref{prop:fixedpoint} to our dual representation
\eqref{eq:mpdiamond}, we see that orbitwise perturbing the right-hand side vector
according $H\leq G$ will preserve the symmetries of $H$.
The perturbed system thus has the form 
\begin{equation}
  \label{eq:polar-perturb}
  Vx\geq -(u + \mu),
\end{equation}
where $\mu_i=\sigma_j \varepsilon_j$ for $\sigma_j\in \{ \pm1 \}$ and 
$j$ is the index of the orbit containing $v_i$.  To ensure an
orbitwise lexicographic perturbation, we will insist
\begin{equation}
  \label{eq:epscond}
  1\gg \varepsilon_1 \gg \varepsilon_2 \gg\dots \gg \varepsilon_k > 0,
\end{equation}
where by $x \gg y$ we mean that $y$ is much smaller than $x$,
i.e.\ it is not possible to combinatorially change the polyhedron defined by 
\eqref{eq:polar-perturb} by choosing $y>0$ smaller.
To implement this symbolically, we need a modification of the
standard lexicographic pivot rule (see~\cite{avis-2000}
or~\cite{chvatal} for more details). Let $b=-(u+\mu)$ and $A=[V \ - \!\! I]$. 
After adding {\em slack variables} to \eqref{eq:polar-perturb}, 
we are left with a system of the form
\begin{equation*}
  Ax= b
\end{equation*}
with $n$ rows and $n+d$ columns, where the first $d$ columns are the
{\em decision variables}.

A \emph{feasible basis} consists of a partition $(\beta,\eta)$ 
of the column indices such that $A_\beta$ (columns of $A$ indexed by $\beta$)
is non-singular and any slack variables in $x^*_\beta=A_\beta^{-1}b$ are non-negative.  
In order to move from one feasible basis to another, we need to perform a \emph{pivot}.  
We start by choosing a column index $j$ to leave $\beta$. 
To find a column index to replace $j$, we need to find
\begin{equation*}
  \argmin_i \frac{b^*_i}{a_i}
\end{equation*}
where $b^*=A^{-1}_\beta b$ and $a=A^{-1}_\beta A_\eta$.  In our case,
where $b=-u-\mu$, we may decompose $b^*$, and thus the ratio test into
two parts $A^{-1}_\beta (-u)$ and $A^{-1}_\beta (-\mu)$.  Because the
values of the $\varepsilon_j$ are chosen very small, the second
part is considered only to break ties.
We write
\begin{equation*}
A^{-1}_\beta (-\mu) = N
\begin{bmatrix}
  \varepsilon_1\\
  \varepsilon_2\\
  \vdots\\
  \varepsilon_k
\end{bmatrix}
\end{equation*}
where column $j$ of $N$ is defined by summing the columns of $A^{-1}_\beta$
corresponding to orbit $j$ of generators, and multiplying by $-\sigma_j$.
Because of the ordering \eqref{eq:epscond},
in order to evaluate 
\begin{equation*}
  \argmin_i \frac{(N\varepsilon)_i}{a_i}
\end{equation*}
we proceed column by column in $N$, reducing the set of ties at each
iteration. 

\subsection{Other refinements and implementation details}
\subsubsection{Adjacency Decomposition Pruning}

Consider facets $F_0$ and $F_1$ that are equivalent under some
symmetry of the basis automorphism group.  This same symmetry acts as
an isomorphism between the corresponding basis graphs. It follows
that when we discover a basis $B$ defining a new orbit, but the facet
$F$ spanned by $B$ is known, we do not need to explore the neighbours
of $B$ since they will be explored in our canonical (i.e.\ discovered
first) facet in the orbit of $F$. In order to ensure orbits are not
discarded, we are careful not to mark $B$ as known until its canonical
discovery.  

\begin{figure}
  \centering
  \includegraphics[height=1.5in]{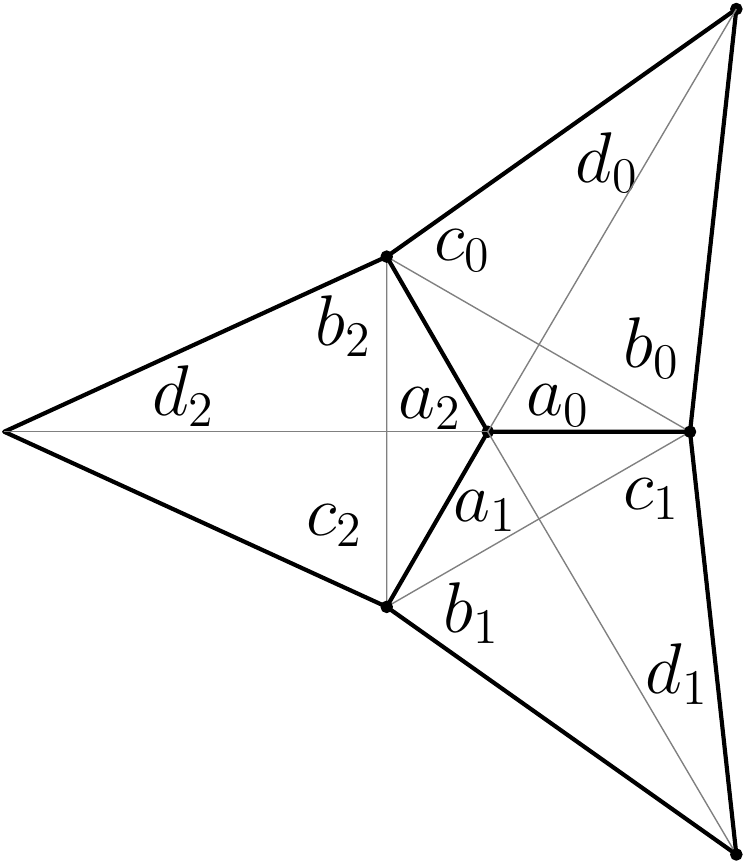}
  \caption{$3$ facets with a $3$-fold rotational symmetry.  Basis
    $q_i$ has vertices $p_i, q_i, r_i$ around the boundary of the
    facet.}
  \label{fig:rotate}
\end{figure}

\begin{figure}
  \centering
  \begin{tabular}{c@{\hspace*{0in}}c}
    \begin{tikzpicture}[anchor=base,>=latex,join=bevel,]
  \pgfsetlinewidth{1bp}
\begin{scope}
  \pgfsetstrokecolor{black}
  \pgfsetstrokecolor{white}
  \pgfsetfillcolor{white}
  \filldraw (0bp,0bp) -- (0bp,126bp) -- (161bp,126bp) -- (161bp,0bp) -- cycle;
\end{scope}
  \draw [->] (95bp,112bp) .. controls (92bp,110bp) and (88bp,107bp)  .. (76bp,98bp);
  \draw [->] (59bp,84bp) .. controls (56bp,82bp) and (54bp,79bp)  .. (43bp,70bp);
  \draw [->] (32bp,56bp) .. controls (31bp,55bp) and (31bp,53bp)  .. (25bp,42bp);
  \draw [->,dashed] (19bp,28bp) .. controls (18bp,27bp) and (18bp,25bp)  .. (14bp,14bp);
  \draw [->,dashed] (25bp,28bp) .. controls (26bp,27bp) and (26bp,25bp)  .. (32bp,14bp);
  \draw [->,dashed] (35bp,56bp) .. controls (35bp,48bp) and (35bp,35bp)  .. (35bp,14bp);
  \draw [->,dashed] (38bp,56bp) .. controls (39bp,55bp) and (39bp,53bp)  .. (45bp,42bp);
  \draw [->,dashed] (26bp,56bp) .. controls (22bp,52bp) and (17bp,47bp)  .. (14bp,42bp) .. controls (12bp,37bp) and (11bp,30bp)  .. (10bp,14bp);
  \draw [->,dashed] (66bp,84bp) .. controls (65bp,69bp) and (61bp,37bp)  .. (56bp,28bp) .. controls (55bp,26bp) and (53bp,24bp)  .. (45bp,14bp);
  \draw [->,dashed] (65bp,84bp) .. controls (62bp,76bp) and (57bp,62bp)  .. (50bp,42bp);
  \draw [->,dashed] (69bp,84bp) .. controls (69bp,83bp) and (70bp,81bp)  .. (72bp,70bp);
  \draw [->,dashed] (75bp,84bp) .. controls (77bp,82bp) and (79bp,80bp)  .. (89bp,70bp);
  \draw [->,dashed] (104bp,112bp) .. controls (104bp,111bp) and (104bp,110bp)  .. (104bp,98bp);
  \draw [->,dashed] (110bp,112bp) .. controls (112bp,110bp) and (113bp,108bp)  .. (121bp,98bp);
  \draw [->,dashed] (115bp,112bp) .. controls (119bp,110bp) and (125bp,106bp)  .. (139bp,98bp);
  \draw [->,dashed] (100bp,112bp) .. controls (99bp,108bp) and (97bp,103bp)  .. (96bp,98bp) .. controls (95bp,92bp) and (95bp,86bp)  .. (96bp,70bp);
\begin{scope}
  \pgfsetstrokecolor{black}
  \draw (104bp,115bp) node {$b_0$};
\end{scope}
\begin{scope}
  \pgfsetstrokecolor{black}
  \draw (67bp,87bp) node {$a_1$};
\end{scope}
\begin{scope}
  \pgfsetstrokecolor{black}
  \draw (35bp,59bp) node {$c_2$};
\end{scope}
\begin{scope}
  \pgfsetstrokecolor{black}
  \draw (22bp,31bp) node {$d_2$};
\end{scope}
\begin{scope}
  \pgfsetstrokecolor{black}
  \draw (11bp,3bp) node {$b_2$};
\end{scope}
\begin{scope}
  \pgfsetstrokecolor{black}
  \draw (35bp,3bp) node {$a_2$};
\end{scope}
\begin{scope}
  \pgfsetstrokecolor{black}
  \draw (48bp,31bp) node {$b_1$};
\end{scope}
\begin{scope}
  \pgfsetstrokecolor{black}
  \draw (74bp,59bp) node {$d_1$};
\end{scope}
\begin{scope}
  \pgfsetstrokecolor{black}
  \draw (97bp,59bp) node {$c_1$};
\end{scope}
\begin{scope}
  \pgfsetstrokecolor{black}
  \draw (104bp,87bp) node {$c_0$};
\end{scope}
\begin{scope}
  \pgfsetstrokecolor{black}
  \draw (127bp,87bp) node {$d_0$};
\end{scope}
\begin{scope}
  \pgfsetstrokecolor{black}
  \draw (150bp,87bp) node {$a_0$};
\end{scope}
\end{tikzpicture} &
   \begin{tikzpicture}[anchor=base,>=latex,join=bevel,]
  \pgfsetlinewidth{1bp}
\begin{scope}
  \pgfsetstrokecolor{black}
  \pgfsetstrokecolor{white}
  \pgfsetfillcolor{white}
  \filldraw (0bp,0bp) -- (0bp,126bp) -- (115bp,126bp) -- (115bp,0bp) -- cycle;
\end{scope}
  \draw [->] (65bp,112bp) .. controls (62bp,110bp) and (60bp,107bp)  .. (49bp,98bp);
  \draw [->] (45bp,84bp) .. controls (46bp,83bp) and (47bp,81bp)  .. (53bp,70bp);
  \draw [->] (57bp,56bp) .. controls (57bp,55bp) and (57bp,54bp)  .. (57bp,42bp);
  \draw [->,dashed] (73bp,112bp) .. controls (73bp,95bp) and (72bp,49bp)  .. (72bp,14bp);
  \draw [->,dashed] (39bp,84bp) .. controls (35bp,71bp) and (25bp,42bp)  .. (16bp,14bp);
  \draw [->,dashed] (41bp,84bp) .. controls (42bp,71bp) and (43bp,43bp)  .. (43bp,14bp);
  \draw [->,dashed] (53bp,28bp) .. controls (52bp,27bp) and (52bp,25bp)  .. (46bp,14bp);
  \draw [->,dashed] (46bp,28bp) .. controls (42bp,26bp) and (38bp,23bp)  .. (25bp,14bp);
  \draw [->,dashed] (61bp,28bp) .. controls (62bp,27bp) and (62bp,25bp)  .. (68bp,14bp);
  \draw [->,dashed] (68bp,28bp) .. controls (72bp,26bp) and (77bp,22bp)  .. (90bp,14bp);
  \draw [->,dashed] (62bp,42bp) .. controls (64bp,46bp) and (67bp,51bp)  .. (68bp,56bp) .. controls (72bp,71bp) and (73bp,89bp)  .. (73bp,112bp);
  \draw [->,dashed] (52bp,42bp) .. controls (50bp,46bp) and (47bp,51bp)  .. (46bp,56bp) .. controls (44bp,61bp) and (43bp,68bp)  .. (41bp,84bp);
  \draw [->,dashed] (59bp,70bp) .. controls (61bp,78bp) and (65bp,92bp)  .. (71bp,112bp);
  \draw [->,dashed] (75bp,112bp) .. controls (79bp,95bp) and (91bp,49bp)  .. (99bp,14bp);
\begin{scope}
  \pgfsetstrokecolor{black}
  \draw (73bp,115bp) node {$b_0$};
\end{scope}
\begin{scope}
  \pgfsetstrokecolor{black}
  \draw (41bp,87bp) node {$c_0$};
\end{scope}
\begin{scope}
  \pgfsetstrokecolor{black}
  \draw (57bp,59bp) node {$d_0$};
\end{scope}
\begin{scope}
  \pgfsetstrokecolor{black}
  \draw (57bp,31bp) node {$a_0$};
\end{scope}
\begin{scope}
  \pgfsetstrokecolor{black}
  \draw (72bp,3bp) node {$a_1$};
\end{scope}
\begin{scope}
  \pgfsetstrokecolor{black}
  \draw (14bp,3bp) node {$a_2$};
\end{scope}
\begin{scope}
  \pgfsetstrokecolor{black}
  \draw (43bp,3bp) node {$b_2$};
\end{scope}
\begin{scope}
  \pgfsetstrokecolor{black}
  \draw (101bp,3bp) node {$c_1$};
\end{scope}
\end{tikzpicture}
  \end{tabular}
  \caption{Comparing the search with pruning (right) and without (left).}
  \label{fig:prune}
\end{figure}
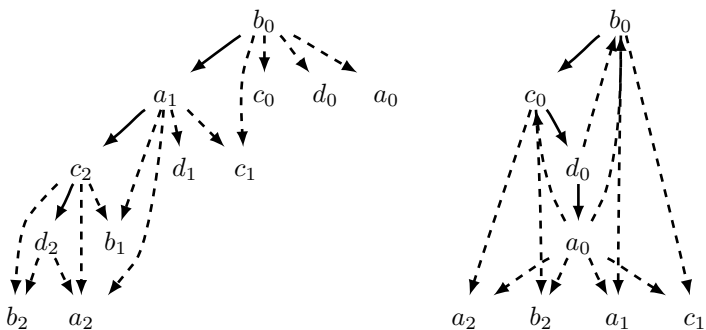

Although this pruning does not reduce the number of orbits of bases
explored, it can reduce the number of actual bases visited (and tested
for isomorphism), since bases of a given orbit are not revisited 
in every copy of the facet $F$. As an example, consider the $3$
quadrilateral facets illustrated in Figure~\ref{fig:rotate}, with
rotational symmetry yielding $4$ orbits of bases.
Without pruning, a depth first search visits all of the bases; with 
pruning only $8$ of $12$ bases are visited (see Figure~\ref{fig:prune}).
The speedup obtained depends roughly on the number of facets visited 
by the unpruned search, which is bounded by the number of basis orbits.

\subsubsection{Metric Invariants}
\label{sec:metric-invariants}
In the case where our symmetry group preserves the inner product
between pairs of vectors, as is the case for the restricted
automorphisms discussed in this paper, we may take advantage of this
in several ways. 

For any face or basis $X$ to be tested for isomorphism, we may
construct a graph (analogous to that constructed in Proposition
\ref{prop:isomorphism-equivalence}) whose nodes are the vectors of $X$
and whose edges are the angles between them.  This graph contains
geometric information not present in the index sets representing $X$,
which can help to speed up an
algorithm to find an isomorphism.  
A simpler observation, and equally widely applicable, is
that the set of pairwise inner products of two isomorphic faces or
bases must be equal. This allows us to store orbit representatives in
a data structure such as a hash table or a balanced tree, with the key
to the data structure being the set of inner products.  This permits
more efficient isometry testing by retrieving exactly those orbit
representatives which pass the inner product invariant.

It is computationally easy to test whether a given linear
transformation $T$ is in the restricted automorphism group $\Aut(V)$ of
vector family $V$.  Since we are interested in restricted
automorphisms carrying basis $X$ to basis Y, we can additionally test
if $\Pi X T = Y$ for some $\Pi \in \Aut(X)$ (where $\Aut(X)$ can be
computed by the same techniques as
Proposition~\ref{prop:isomorphism-equivalence}). We have only
implemented an exhaustive search of $\Aut(X)$, and this is naturally
only effective when $\Aut(X)$ is quite small.  In principle it should
be possible to integrate the test for $T$ being a restricted
isomorphism into a backtracking procedure to search for $\Pi$.

\section{Conclusions}

Much as in the case of polyhedral representation conversion without
symmetries, a certain amount of trial and error seems to be necessary
to decide on the the best method to attack a given conversion problem
up to symmetries. Currently decomposition methods have the best record
of solving interesting problems; on the other hand current software
requires a certain amount of user intervention in the form of choosing
how to treat subproblems. It would be helpful to automate this
process. In this context, a virtue of the pivoting methods is that
good methods to estimate their running time exist~\cite{ad94}. It
would be beneficial, not just when working with symmetry, to have
effective methods (or at least heuristics) for estimating the running
time of incremental methods.

\section*{Acknowledgements}

The authors would like to acknowledge fruitful discussions with David
Avis, Antoine Deza, Komei Fukuda, Alexander Hulpke, Michael Joswig,
Jesus de Loera, Brendan McKay, Hugh Thomas, and Frank Vallentin.  They
would also like to thank the Centre de Recherche Math\'ematiques for
making possible the workshop at which many of these discussions took
place.

\end{document}